\RequirePackage{fix-cm}
\RequirePackage{luatex85}
\documentclass[3p]{elsarticle}
\usepackage[utf8]{inputenc}
\usepackage[T1]{fontenc}
\usepackage{lmodern}
\usepackage{mathtools}
\usepackage{amsfonts}
\usepackage{amssymb}
\usepackage{amsthm}
\usepackage{color}
\usepackage{multirow}
\usepackage{array}
\usepackage{booktabs}
\usepackage{faktor}
\usepackage{cases}
\usepackage{fontawesome5}
\usepackage[final]{microtype}
\usepackage{bm}

\usepackage{pgfplots} 
\pgfplotsset{compat=1.18}

\usepackage{pgfplotstable}
\usepackage{subcaption}

\definecolor{graph_1}{RGB}{117,112,179}
\definecolor{graph_2}{RGB}{217,95,2}
\definecolor{graph_3}{RGB}{27,158,119}
\definecolor{graph_4}{RGB}{231,41,138}
\definecolor{graph_5}{RGB}{20,170,240}
\definecolor{graph_6}{RGB}{165,42,42}

\newcommand{\impCFL}{{\mu_\varepsilon}}
\newcommand{\lp}{\left(}
\newcommand{\rp}{\right)}
\newcommand{\sfr}[2]{{\scriptstyle \faktor {#1} {#2}}}        

 \newtheorem{algorithm}{Algorithm}

\usepackage{hyperref}
\hypersetup{
    colorlinks=true,%
    linkcolor=blue,%
    citecolor=blue,%
    filecolor=blue,%
    urlcolor=blue,%
    pdfcreator=LaTeX,%
    breaklinks=true,%
    pdfpagelayout=SinglePage,%
    bookmarksopen=true,%
    bookmarksopenlevel=2,
    hypertexnames=true
}
\usepackage{bookmark}

\makeatletter

\newcommand{\printslope}[4]{
   \tikzset{fixed point arithmetic}
   \def\nero@printslope@orderlist{#1}
   \edef\nero@printslope@xpos{#2}
   \edef\nero@printslope@ypos{#3}
   \edef\nero@printslope@width{#4}
   \pgfmathparse{\nero@printslope@xpos+\nero@printslope@width}
   \edef\nero@printslope@px{\pgfmathresult}
   \edef\nero@printslope@py{\nero@printslope@ypos}
   \edef\nero@printslope@qx{\nero@printslope@xpos}
   \edef\nero@printslope@ry{\nero@printslope@ypos}
   \foreach \nero@printslope@order in {#1}{
      \pgfmathparse{
      ((\nero@printslope@px/\nero@printslope@xpos)^(\nero@printslope@order))*\nero@printslope@ypos}
      \edef\nero@printslope@qy{\pgfmathresult}
      \edef\nero@aux1{\noexpand\draw[line width=0.6pt]
         (axis cs:\nero@printslope@xpos,\nero@printslope@ry)
         -- (axis cs:\nero@printslope@qx,\nero@printslope@qy)
         -- (axis cs:\nero@printslope@px,\nero@printslope@py);}
      \nero@aux1
      \pgfmathparse{10^((ln(\nero@printslope@ry)+ln(\nero@printslope@qy))/(ln(10)*2))}
      \edef\nero@printslope@labelpos{\pgfmathresult}
      \edef\nero@aux2{\noexpand\node[anchor=east] at
         (axis cs:\nero@printslope@qx,\nero@printslope@labelpos)
         {\noexpand\small \nero@printslope@order};}
      \nero@aux2
      \global\edef\nero@printslope@ry{\nero@printslope@qy}
   }
   \draw[line width=0.6pt] (axis cs:\nero@printslope@xpos,\nero@printslope@ypos)
      |- (axis cs:\nero@printslope@px,\nero@printslope@py);
   \pgfmathparse{10^((ln(\nero@printslope@px)+ln(\nero@printslope@xpos))/(ln(10)*2))}
   \edef\nero@printslope@labelpos{\pgfmathresult}
   \node[anchor=north] at (axis cs:\nero@printslope@labelpos,\nero@printslope@ypos) {\small 1};
}

\newcommand{\printslopeintwodimensions}[4]{
   \def\nero@printslope@orderlist{#1}
   \edef\nero@printslope@xpos{#2}
   \edef\nero@printslope@ypos{#3}
   \edef\nero@printslope@width{#4}
   \pgfmathparse{\nero@printslope@xpos+\nero@printslope@width}
   \edef\nero@printslope@px{\pgfmathresult}
   \edef\nero@printslope@py{\nero@printslope@ypos}
   \edef\nero@printslope@qx{\nero@printslope@xpos}
   \edef\nero@printslope@ry{\nero@printslope@ypos}
   \foreach \nero@printslope@order in {#1}{
      \edef\nero@printslope@qy{\fpeval{((\nero@printslope@px/\nero@printslope@xpos)^(\nero@printslope@order/2))*\nero@printslope@ypos}}
      \edef\nero@aux1{\noexpand\draw[line width=0.7pt]
         (axis cs:\nero@printslope@xpos,\nero@printslope@ry)
         -- (axis cs:\nero@printslope@qx,\nero@printslope@qy)
         -- (axis cs:\nero@printslope@px,\nero@printslope@py);}
      \nero@aux1
      \edef\nero@printslope@labelpos{\fpeval{10^((ln(\nero@printslope@ry)+ln(\nero@printslope@qy))/(ln(10)*2))}}
      \edef\nero@aux2{\noexpand\node[anchor=east] at
         (axis cs:\nero@printslope@qx,\nero@printslope@labelpos)
         {\noexpand\small \nero@printslope@order};}
      \nero@aux2
      \global\edef\nero@printslope@ry{\nero@printslope@qy}
   }
   \draw[line width=0.7pt] (axis cs:\nero@printslope@xpos,\nero@printslope@ypos)
      |- (axis cs:\nero@printslope@px,\nero@printslope@py);
   \pgfmathparse{10^((ln(\nero@printslope@px)+ln(\nero@printslope@xpos))/(ln(10)*2))}
   \edef\nero@printslope@labelpos{\pgfmathresult}
   \node[anchor=north] at (axis cs:\nero@printslope@labelpos,\nero@printslope@ypos) {\small 1};
}

\makeatother

\newtheorem{theorem}{Theorem}
\newtheorem{lemma}[theorem]{Lemma}

\theoremstyle{remark}
\newtheorem{remark}{Remark}

\theoremstyle{definition}

\journal{Applied Mathematics and Computation}

\title{TVD-MOOD schemes based on implicit-explicit time integration}

\date{July 4, 2022}

\author[1]{Victor Michel-Dansac}
\ead{victor.michel-dansac@inria.fr}
\author[2]{\texorpdfstring{Andrea Thomann \corref{cor1}}{Andrea Thomann}}
\ead{athomann@uni-mainz.de}

\address[1]{Université de Strasbourg, CNRS, Inria, IRMA, F-67000 Strasbourg, France}
\address[2]{Institut f\"ur Mathematik, Johannes Gutenberg-Universit\"at Mainz, Germany}

\cortext[cor1]{Corresponding author}

\numberwithin{equation}{section}

\begin{document}

\newlength{\plotwidth}
\setlength{\plotwidth}{0.95\linewidth/2}
\newlength{\plotheight}
\setlength{\plotheight}{0.7\linewidth/2}
\newlength{\plotlinewidth}
\setlength{\plotlinewidth}{1.2pt}
\newlength{\smallerplotlinewidth}
\setlength{\smallerplotlinewidth}{0.8pt}

\begin{abstract}
	The context of this work is the development of first order total variation diminishing (TVD) implicit-explicit (IMEX) Runge-Kutta (RK) schemes as a basis of a Multidimensional Optimal Order detection (MOOD) approach to approximate the solution of hyperbolic multi-scale equations.
	A key feature of our newly proposed TVD schemes is that the resulting CFL condition does not depend on the fast waves of the considered model, as long as they are integrated implicitly.
	However, a result from Gottlieb et al.~\cite{GotShuTad2001} gives a first order barrier for unconditionally stable implicit TVD-RK schemes and TVD-IMEX-RK schemes with scale-independent CFL conditions.
	Therefore, the goal of this work is to consistently improve the resolution of a first-order IMEX-RK scheme, while retaining its $L^\infty$ stability and TVD properties.
	In this work we present a novel approach based on a convex combination between a first-order TVD IMEX Euler scheme and a potentially oscillatory high-order IMEX-RK scheme.
	We derive and analyse the TVD property for a scalar multi-scale equation and
	numerically assess the performance of our TVD schemes compared to standard $L$-stable and SSP IMEX RK schemes from the literature.
	Finally, the resulting TVD-MOOD schemes are applied to the isentropic Euler equations.
\end{abstract}


\maketitle

\section{Introduction}
Multi-scale equations arise in a wide range of applications,
such as shallow flows \cite{BisAruLukNoe2014},
magnetohydrodynamics~\cite{MatBro1988}, multi-material interfaces~\cite{AbbIollPupp2019} or atmospheric flows~\cite{Kle2010}.
When developing numerical methods for such applications,
it is of prime importance to obtain physically admissible solutions under these multi-scale constraints.
%
%
In order to numerically treat these different scales, one must assess whether the fast scales are relevant to the physical solution.
Indeed, accurately capturing these fast scales requires a very restrictive time step.
This issue is discussed e.g. in~\cite{GuiVio999} for the Euler equations.
When the impact of the fast scales on the physical solution is less important,
numerical methods that do not accurately capture all scales, but only follow the slow dynamics, are necessary.
One option, which we will study in this paper, is to use Implicit-Explicit (IMEX) schemes, where the terms associated to the fast wave propagation are treated implicitly.
Such schemes are well-studied in the literature,
see for instance~\cite{AscRuuSpi1997} for efficient IMEX schemes applied to hyperbolic-parabolic problems,~\cite{ParRus2005} for IMEX schemes adapted to stiff relaxation source terms,
or~\cite{NoeBisAruLukMun2014,DimLouVig2017,BosRusSca2018,ParMun2005,BusRioMVazCDum2021}, and references therein, for IMEX schemes designed for the low Mach regime of the Euler equations.
In this work, we are concerned with hyperbolic systems whose stiffness originates from the flux, rather than a source term.
Let us emphasise that we will not consider hyperbolic systems with stiff source terms typically arising from relaxation processes.
For their treatment, we refer for instance to~\cite{ParRus2005}.

%
To increase the quality of numerical approximations,
one may turn to high-order schemes.
However, such schemes are known to introduce spurious oscillations in the solution away from smooth regions.
This is problematic, especially when considering non-linear hyperbolic equations,
as the solution can develop discontinuities even when starting with a smooth initial condition.
This was already observed by Harten in~\cite{Har1984}, who introduced the notion of total variation diminishing (TVD) schemes,
and constructed non-oscillatory explicit and implicit second-order TVD schemes.
Those schemes are non-linear, even when applied on linear equations,
as from Godunov's theorem~\cite{God1959} it follows that linear TVD schemes can only be first-order accurate.
Since non-linear implicit schemes are computationally very costly, especially when applied to non-linear systems of equations,
the construction of higher order explicit TVD schemes
remained an active area of research, see e.g.~\cite{Swe1984,ShuOsh1988,GotShu1998,BreGotGraHigKetNe2017} and references therein.
Later, in the more general framework of strong stability preserving (SSP) implicit and explicit schemes~\cite{GotShuTad2001},
the stability property is achieved by relying on convexity arguments regarding forward and backward Euler schemes,
rather than adding artificial viscosity to achieve the TVD property, as was done in~\cite{Har1984,Swe1984,Roe1984}.
The high-order explicit and implicit SSP schemes developed in~\cite{GotShuTad2001,Got2005,GotKetShu2011}
are limited by a very restrictive CFL condition comparable to a forward Euler scheme in order to remain oscillation-free.
This makes the use of high-order implicit SSP schemes rather costly and impractical in applications,
compared to high-order explicit SSP schemes, as was remarked in~\cite{GotKetShu2011}.
Regarding IMEX SSP schemes, we refer for instance to~\cite{HigHapKocKup2014,ConGotGraSha2017,HigKetKoc2018}.
All high-order SSP schemes mentioned above require the time step to depend on all scales to achieve stability, but are provably high-order accurate.
Unfortunately, they are not well-suited to the multi-scale setting,
where the time step is strongly restricted by the fast scale leading, in extreme cases, to a vanishing time step.
We also want to mention the very recent work~\cite{GotGraHuShu2022},
whose authors derive a high-order IMEX SSP scheme with a scale-independent time step restriction.
However, this scheme cannot be applied unless some restrictive assumptions are satisfied by the system under consideration.

In contrast, our focus here is the construction of first order IMEX TVD schemes,
whose CFL restriction solely stems from the explicitly treated terms associated to a scale-independent material velocity.
The work presented in this manuscript is greatly motivated by the seminal work by Gottlieb et al.~\cite{GotShuTad2001},
where it was proven that an unconditionally TVD implicit RK scheme is at most first-order accurate,
see also \cite{Spi1983,Hig2006}.
Unfortunately, this result holds also for IMEX discretisations with a scale-independent CFL restriction.
In fact, this discouraging result is also observed in~\cite{DimLouMicVig2018,BouFraNav2020}
when attempting to construct second-order TVD IMEX schemes for the Euler equations.

Although the TVD property is crucial to accurately capture discontinuities in the numerical solution, it becomes of less importance in smooth regions.
In order to achieve a high-order approximation of the solution in such regions,
while keeping the solution as oscillation-free as possible in the vicinity of discontinuities,
we adapt the IMEX framework to a procedure inspired from the MOOD
(Multidimensional Optimal Order Detection)
techniques from~\cite{ClaDioLou2011}.
The gist of the MOOD framework is to lower
the order of accuracy of the scheme near problematic zones,
i.e. areas where the high-order scheme violates some predetermined admissibility constraint.
Therefore, a lower order scheme with good stability properties,
called a \emph{parachute scheme}, is needed in the MOOD framework.

In the present work, as stated before, we design first-order TVD IMEX RK discretisations that are consistently less diffusive than the standard first-order backward/forward Euler (IMEX1) scheme.
Such discretisations therefore provide suitable fall-back schemes for the MOOD approach,
yielding a reduced space-time error compared to using a traditional IMEX1 scheme as a parachute.
The approach given here to construct such a fall-back scheme builds on the results from~\cite{DimLouMicVig2018,MicTho2019},
where the increase in precision was achieved by introducing a convex combination of said first-order TVD scheme with an oscillatory second-order scheme.
In~\cite{DimLouMicVig2018}, the ARS(2,2,2) scheme from~\cite{AscRuuSpi1997} was used as a basis for the convex combination,
and this result was extended to a general class of second-order IMEX RK schemes in~\cite{MicTho2019}.
Here, we generalize and extend the results from~\cite{DimLouMicVig2018,MicTho2019} even further,
applying a convex combination to each stage,
rather than merely to the final update as done in \cite{DimLouMicVig2018,MicTho2019}.
This allows to construct TVD schemes based on arbitrarily high order IMEX RK schemes.
Note that convex combinations have already been used to recover first-order properties lost at higher orders,
see for instance~\cite{HuAdaShu2013} to recover the positivity property or~\cite{MicBerClaFou2017} for well-balanced problems.

%
%
The paper is organised as follows.
In Section~\ref{sec:problem_description},
we motivate the problem of multi-scale equations, illustrated by a scalar linear transport equation.
We also introduce notation for the space discretisation.
Section \ref{sec:scheme_derivation} revisits the MOOD procedure and details the construction of the numerical scheme.
The formalism of IMEX-RK is briefly recalled and the TVD constraints are derived for the time-semi discrete scheme.
Subsequently, the fully discrete scheme is given, based on a finite volume approach and TVD limiters.
To completely determine the scheme, the choice of free parameters in the convex combinations is addressed.
Section~\ref{sec:numerics} is devoted to numerical experiments to verify the properties of the parachute first order TVD scheme as well as the MOOD scheme.
To numerically validate that our TVD-IMEX-MOOD schemes are a noticeable improvement over widely used $L$-stable IMEX and the SSP IMEX schemes,
we compare the performance of the schemes in terms of accuracy, CPU times and CFL restrictions on discontinuous solutions of the scalar multi-scale equation.
Moreover, we numerically show that using our TVD scheme as a basis of the MOOD procedure shows a significant reduction of the space-time error.
We finally apply the scheme to the isentropic Euler equations, assessing its performance for Riemann problems in two space dimensions and an explosion problem,
as well as verifying its accuracy and the convergence in the singular Mach number limit.
To complete this manuscript, a conclusion is presented in Section~\ref{sec:conclusion}.


\section{Motivation}
\label{sec:problem_description}

To study conditions to obtain a TVD scheme, we consider the linear advection equation
\begin{equation}
	\label{eq:ScalarMultiscale}
	\begin{dcases}
		w_t + c_m w_x + \frac{c_a}{\varepsilon} w_x = 0, \\
		w(0,x) = w^0(x),
	\end{dcases}
\end{equation}
where $w: (\mathbb{R}^+, \Omega) \to \mathbb{R}$, $\Omega \subset \mathbb{R}$.
In \eqref{eq:ScalarMultiscale}, $c_m$ and $c_a / \varepsilon$ denote two transport speeds which can differ significantly depending on the choice of the parameter $\varepsilon > 0$.
Without loss of generality, we consider only positive transport directions, i.e. $c_m > 0$ and $c_a > 0$.
In the following, we consider $c_m = \mathcal{O}(1)$ and $0 < c_a/\varepsilon \leq c_m$.
The term $\frac{c_a}{\varepsilon} w_x$ is stiff as soon as $\varepsilon \ll 1$,
and applying a purely explicit scheme would lead to a severe time step restriction, given by
\begin{equation}
	\label{eq:CFLac}
	\Delta t \leq \varepsilon \ \nu_{\text{ac}} \ \frac{\Delta x}{\varepsilon c_m + c_a},
\end{equation}
with a CFL coefficient $\nu_{\text{ac}}$ independent of $\varepsilon$.
However, when $\varepsilon$ tends to zero, the above time step vanishes,
leading to huge computational costs, especially when considering long time periods.
Therefore, to avoid such a restriction,
we integrate the wave associated to $c_a/\varepsilon$ implicitly,
whereas $c_m w_x $ is treated explicitly.
This approach leads to a CFL condition oriented to the slow wave $c_m$,
independently of $\varepsilon$.
It is given by
\begin{equation}
	\label{eq:CFLmat}
	\Delta t \leq \nu_{\text{mat}} \ \frac{\Delta x}{c_m}.
\end{equation}
In space, we use an upwind discretization since it was shown in \cite{DimLouVig2017} that the use of centred differences destroys the $L^\infty$ stability for non-linear systems.
Even though our considered problem is linear, the goal is to apply the scheme on non-linear systems such as the isentropic Euler equations discussed in Section \ref{sec:isentropic_Euler}.

\begin{figure}[!tb]

	\centering
	\makebox[0.95\linewidth][c]{

		\includegraphics{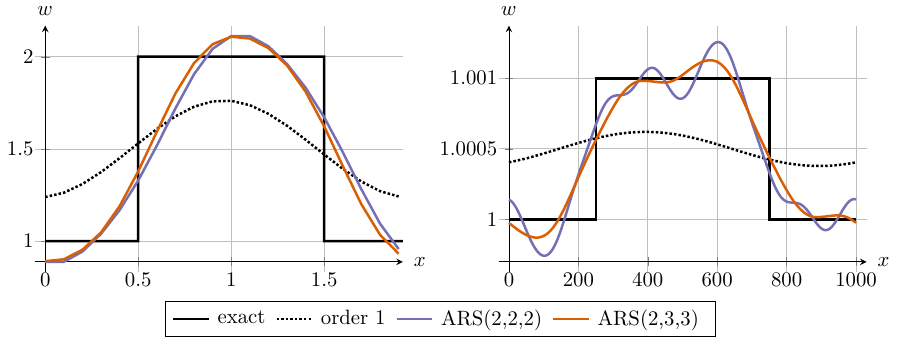}

	}

	\caption{
		Approximation of a discontinuous solution with $\Delta x = 0.1$
		using the first-order, second-order ARS(2,2,2) and third-order ARS(2,3,3) scheme for $\varepsilon = 1$ (left) and $\varepsilon = 10^{-3}$ (right),
		with an upwind space discretisation.
	}

	\label{fig:motivation}

\end{figure}

The space and time discretisations follow the usual finite volume framework,
where the computational domain~$\Omega$ is divided in $N$ uniformly spaced cells $C_j = (x_{j-1/2},x_{j+1/2})$, of size $\Delta x$ and whose center is $x_j = j \Delta x$.
The solution in cell $C_j$ is then approximated by the cell average, given by
\begin{equation}
	w_j(t) \approx \frac{1}{\Delta x}\int_{\Omega_j} w(x,t) dx.
\end{equation}
A first order space semi-discrete scheme
approximating weak solutions to~\eqref{eq:ScalarMultiscale} is then given by
\begin{equation}
	\label{eq:SpaceDiscFirst}
	\partial_t w_j(t) + \frac{c_m}{\Delta x} \Delta_j(t) + \frac{c_a}{\varepsilon \Delta x} \Delta_j(t) = 0,
\end{equation}
where we have introduced the abbreviation $\Delta_j(t) = w_j(t) - w_{j-1}(t)$.
To obtain a fully discrete scheme, a suitable implicit-explicit time integration method
has to be applied.
We discretise the time variable with $t^{n+1} = t^n + \Delta t$,
where $\Delta t$ denotes the time step, which has to obey a material CFL condition \eqref{eq:CFLmat}.
However, it is well-known, see e.g. \cite{GotShuTad2001}, that approximating discontinuous
solutions with high-order non-TVD methods can lead to spurious artefacts near jump positions.
This behaviour is illustrated in Figure~\ref{fig:motivation},
where we display the approximation of an advected rectangular
bump profile with the first-order scheme
\begin{equation}
	\label{eq:firstOrder}
	w^{n+1}_j = w^n_j - c_m \frac{\Delta t}{\Delta x} \Delta_j^n - \frac{c_a}{\varepsilon} \frac{\Delta t}{\Delta x} \Delta_j^{n+1},
\end{equation}
as well as with the well-known second-order ARS(2,2,2) and third-order ARS(2,3,3)
IMEX schemes from~\cite{AscRuuSpi1997}.
The details on the numerical experiment are given in Section~\ref{sec:numerics_toy_problem}.
We clearly observe in Figure~\ref{fig:motivation} that the ARS methods
violate the bounds on the numerical solution in both cases.
Therefore, in order to avoid oscillations as in Figure~\ref{fig:motivation},
we need $L^\infty$ stable or TVD schemes.
A scheme is said to be~$L^\infty$ stable if
\begin{equation}
	\label{def:LinftyStab}
	\| w^{n+1}\|_\infty = \underset{ j \in \lbrace 1, \ldots, N \rbrace }{\max} |w_j^{n+1}| \leq \|w^n\|_\infty,
\end{equation}
and TVD if
\begin{equation}
	\label{def:TVD}
	\text{TV}(w^{n+1}) = \sum_{j=1}^{N} \left|w_{j+1}^{n+1} - w_j^{n+1}\right| \leq \text{TV}(w^n).
\end{equation}

Unfortunately, it can be proven for IMEX RK schemes,
following a result of Gottlieb et al.~\cite{GotShuTad2001},
that there cannot exist $L^\infty$ stable IMEX RK schemes of order $p \geq 2$
whose CFL restriction only stems from the explicitly treated part.
Therefore, we do not look for higher order TVD IMEX integrators.
However, we clearly see in Figure~\ref{fig:motivation} that the first-order scheme is too diffusive to capture the structure of the solution, and thus it is also not a good choice as a base scheme in the MOOD hierarchical structure.
As a consequence, our main focus here is to construct a first-order IMEX integration scheme fulfilling the $L^\infty$ stability property~\eqref{def:LinftyStab} and the TVD property~\eqref{def:TVD},
that has a reduced numerical diffusion compared to the first order scheme \eqref{eq:firstOrder}, and is therefore suited as a fallback scheme when combined with a MOOD procedure.

\section{Derivation of the numerical scheme}
\label{sec:scheme_derivation}

The goal of this section is the derivation of a numerical scheme based on a MOOD-like procedure.
The usual MOOD framework for explicit schemes, see e.g.~\cite{ClaDioLou2011},
consists in locally and gradually lowering the order of the scheme when an oscillation is detected.
The lowest order scheme the procedure can use is called a \textit{parachute scheme}.
Generally speaking, it guarantees the preservation of desired properties not satisfied by the higher order schemes.
Here, the higher order scheme is given by a standard IMEX-RK scheme,
and the lowest order scheme consists of a first order TVD-IMEX scheme.
Consequently, the precision of the first order TVD scheme can be improved without degrading its stability properties, yielding a high order approximation in smooth flow regimes.
Due to the non-local nature of the implicit part in the IMEX schemes, the standard MOOD algorithm as given in \cite{ClaDioLou2011} has to be modified, since the solution cannot be recomputed on a few selected cells only.
Therefore, it has to be updated on the whole mesh.

Depending on the application, different detection criteria can be applied to identify oscillations violating the TVD property.
For the toy problem \eqref{eq:ScalarMultiscale},
an oscillation is detected when the solution leaves the bounds of the initial condition.
Indeed, the unknown $w$ for this toy problem satisfies
the maximum principle $\partial_t \| w(t, x) \|_\infty \leq 0$,
and therefore an $L^\infty$-stable discretisation obeys
$\| w^{n+1} \|_\infty \leq \| w^0 \|_\infty$ for all $n \geq 0$.
Thus, the detection criterion, denoted by $\Phi$, depends on the initial condition only and is time independent.
However, in general, time dependent detection criteria must be taken into account for non-linear problems, as detailed for the isentropic Euler equations in Section~\ref{sec:numerics}.
The modified MOOD scheme, see also \cite{DimLouMicVig2018,MicTho2019}, is given by the following algorithm.
\begin{algorithm}[MOOD$p$ scheme]
	\label{algo:implicit_MOOD}
	Define the initial detection criterion given by $\mathcal{E}^0 = \|\Phi (w^0)\|$.
	Equipped with a stable first order TVD scheme, the MOOD$p$ scheme consists in applying the following procedure at each time step:
	\begin{enumerate}
		\item[(1)] Compute a candidate numerical solution $w_c^{n+1}$ with a
			$p^{\text{th}}$ order IMEX-RK scheme.
		\item[(2)] Detect whether an oscillation is present somewhere in the space domain,
			i.e. whether the detection criterion is satisfied by the candidate solution:
			\begin{equation*}
				\tag{\text{DMP}}
				\label{eq:maximum_principle_satisfied}
				\|\Phi( w_c^{n+1}) \|_\infty \leq \mathcal{E}^n.
			\end{equation*}
		\item[(3a)] If~\eqref{eq:maximum_principle_satisfied} holds,
			set for the numerical solution at the new time step $w^{n+1}= w_c^{n+1}$.
		\item[(3b)] Otherwise, compute a solution $w_\text{\faParachuteBox}^{n+1}$ with the first order parachute scheme on the whole computational domain and set $w^{n+1} = w_\text{\faParachuteBox}^{n+1}$.
		\item[(4)]
			For a time-dependent detection criterion,
			update $\mathcal{E}^{n+1} = \xi \| \Phi(w^{n+1}) \|_\infty + (1 - \xi) \mathcal{E}^n$
			with $\xi \in [0,1]$,
			otherwise set $\mathcal{E}^{n+1} = \mathcal{E}^0$.
	\end{enumerate}
\end{algorithm}
Note that, in Step~(4) of the MOOD algorithm,
the detection criterion is relaxed with a convex combination of parameter $\xi \in [0, 1]$
between the current solution and the criterion at the previous time step.
This allows a finer control of the permitted oscillations in the MOOD solution.

A crucial part of the MOOD scheme is the construction of the parachute scheme,
which should preserve the TVD property of the solution
in case the higher order scheme produces oscillations.
This is addressed in the following section.

\subsection{Construction of the time-semi discrete parachute scheme}
\label{sec:parachute_semi_discrete}

The proposed parachute schemes are based on an IMEX-RK approach
where the implicit part is diagonally implicit.
The time update for an~$s$-stage IMEX-RK scheme for equation~\eqref{eq:SpaceDiscFirst} is given by
\begin{equation}
	\label{eq:IMEXupdate}
	w^{n+1}_j = w^n_j - \lambda \sum_{k=1}^{s} \tilde b_k \Delta_j^{(k)} - \impCFL \sum_{k=1}^{s} b_k \Delta_j^{(k)},
\end{equation}
where we have set
\begin{equation*}
	\lambda = \frac{\Delta t}{\Delta x} c_m
	\text{, \qquad}
	\impCFL = \frac{\Delta t}{\Delta x}\frac{c_a}{\varepsilon}
	\text{, \qquad}
	\Delta_j^{(k)} = w_j^{(k)} - w_{j-1}^{(k)},
\end{equation*}
and where the stages are defined as
\begin{equation}
	\label{eq:IMEXstages}
	w^{(k)}_j = w^n_j - \lambda \sum_{l=1}^{k-1} \tilde a_{kl} \Delta_j^{(l)} - \impCFL \sum_{l=1}^{k} a_{kl} \Delta_j^{(l)}.
\end{equation}
For the sake of clarity,
we consider an IMEX-RK method of type CK (Carpenter and Kennedy)~\cite{KenCar2003},
i.e. we take $a_{11} = 0$.
For results based on a non-CK scheme with $a_{11} \neq 0$
see~\ref{sec:appendix_non_CK_schemes}.
In order to obtain a CFL condition like \eqref{eq:CFLmat} which does not depend on $\varepsilon$,
the weights $(a_{k1})_{k \in \left\{ 2, \ldots, s \right\}}$, as well as $b_1$, have to be zero.
This was shown in detail in \cite{MicTho2019} for a generic second-order CK method.
We summarize the structure of the RK scheme in the following Butcher tableaux notation
\begin{equation}
	\label{tab:general}
	\text{explicit: }
	\renewcommand{\arraystretch}{1.25}
	\begin{array}{c|ccccc}
		0          & 0             & 0      & \cdots           & 0            \\
		\tilde c_2 & \tilde a_{21} & 0      & \cdots           & 0            \\
		\vdots     & \vdots        & \ddots & \ddots           & \vdots       \\
		\tilde c_s & \tilde a_{s1} & \cdots & \tilde a_{s,s-1} & 0            \\ \hline
		           & \tilde b_1    & \cdots & \tilde b_{s-1}   & \tilde b_{s}
	\end{array}
	\hskip1cm
	\text{implicit: }
	\renewcommand{\arraystretch}{1.25}
	\begin{array}{c|cccc}
		0      & 0      & 0      & \cdots & 0      \\
		c_2    & 0      & a_{22} & \cdots & 0      \\
		\vdots & \vdots & \vdots & \ddots & \vdots \\
		c_s    & 0      & a_{s2} & \cdots & a_{ss} \\ \hline
		       & 0      & b_{2}  & \cdots & b_{s}
	\end{array},
\end{equation}
with the coefficients $\tilde c$ and $c$ obeying
\begin{equation}
	\tilde c_i = \sum_{j=1}^{i-1} \tilde a_{ij}
	\text{\quad and \quad}
	c_i = \sum_{j=1}^i a_{ij}.
	\label{cond:c}
\end{equation}
Since the TVD scheme we construct are based on higher order IMEX-RK schemes,
the weights have to fulfil high-order compatibility conditions as given in \cite{ParRus2001}.
For orders higher than three, we refer to the order conditions in~\cite{KenCar2003}.

To obtain a TVD scheme from an $s$-stage IMEX-RK scheme,
we propose a convex combination of each stage with
the first-order IMEX Euler scheme~\eqref{eq:firstOrder}
at time $t^n + c_k \Delta t$ for the $k^{\text{th}}$ stage.
This approach yields $s$ free parameters $\theta_k \in [0,1]$,
where $k=1,\ldots, s$ denotes the stage in the IMEX scheme.
The closer each $\theta_k$ is to $1$, the higher the contribution of the IMEX-RK scheme.
As a consequence, we shall seek to maximize the values of $\theta_k$,
to reduce diffusion as much as possible compared to the first-order IMEX Euler scheme.
The stages of the TVD scheme are defined by
\begin{equation}
	\label{eq:IMEXconvexStages}
	w^{(k)}_j
	+
	(1 - \theta_k)c_k \impCFL \Delta_j^{(k)}
	=
	w^n_j
	-
	\lambda \left(
	(1- \theta_k) \tilde c_k \Delta_j^{n}
	+
	\theta_k\sum_{l=1}^{k-1} \tilde a_{kl} \Delta_j^{(l)}
	\right)
	-
	\impCFL \theta_k\sum_{l=1}^{k} a_{kl}\Delta_j^{(l)},
\end{equation}
and the update by
\begin{equation}
	\label{eq:IMEXconvexUpdate}
	w^{n+1}_j
	=
	w^n_j
	- \theta_{s+1} \left(
	\lambda \sum_{k=1}^{s} \tilde b_k \Delta_j^{(k)}
	+
	\impCFL \sum_{k=1}^{s} b_k \Delta_j^{(k)}
	\right)
	-
	( 1 - \theta_{s+1} ) \lp \lambda \Delta_j^n + \impCFL \Delta_j^{n+1} \rp.
\end{equation}
Note that, for Butcher tableaux like~\eqref{tab:general}, we immediately set $\theta_1 = 1$ to recover~$w^{(1)} = w^n$.
This means that the convex stages appear earliest for $k=2$.
In the case where the weights $\tilde b$ and $b$ respectively coincide with the last row of $\tilde A$ and $A$,
the stage $w^{(s)}$ coincides with the final update $w^{n+1}$.
In particular, we then have $\theta_{s+1} = \theta_s$.

We emphasise that a TVD scheme resulting from the
convex combination~\eqref{eq:IMEXconvexUpdate} is at most first-order accurate,
see for instance \cite{Spi1983,GotShuTad2001,Hig2006}.
However, it is a great fit as a parachute scheme for the MOOD Algorithm~\ref{algo:implicit_MOOD}.
Section~\ref{sec:third_order_convex} is dedicated to the necessary conditions
on the free parameters in order to obtain the TVD property using a three stage IMEX-RK scheme.
Then, in Section~\ref{sec:nth_order_convex},
we state an extension to construct TVD discretisations
based on more general Butcher tableaux~\eqref{tab:general}.

\subsubsection{TVD conditions for a three stage scheme}
\label{sec:third_order_convex}

As basis of the TVD scheme,
we consider Butcher tableaux with $s_{\text{eff}} = 3$ effective computational steps.
In the spirit of \eqref{tab:general}, they are given by
\begin{equation}
	\label{tab:ARS4}
	\renewcommand{\arraystretch}{1.25}
	\text{ explicit: }
	\begin{array}{c|cccc}
		0   & 0             & 0             & 0     \\
		c_2 & \tilde a_{21} & 0             & 0     \\
		c_3 & \tilde a_{31} & \tilde a_{32} & 0     \\
		\hline
		    & 0             & b_{2}         & b_{3}
	\end{array}, \hskip0.5cm
	\text{ implicit: }
	\begin{array}{c|cccc}
		0   & 0 & 0      & 0      \\
		c_2 & 0 & a_{22} & 0      \\
		c_3 & 0 & a_{32} & a_{33} \\
		\hline
		    & 0 & b_{2}  & b_{3}
	\end{array}.
\end{equation}
Therein, we have assumed that $\tilde b = b$ and $\tilde c = c$, see also~\cite{ParRus2005}.
This has the advantage that the weights in the final update coincide with respect to the explicitly and implicitly treated terms.
Applying the third-order conditions as given in \cite{ParRus2001} on
the weights and coefficients given by~\eqref{tab:ARS4} lead to the following tableaux,
with $\gamma \notin \{ 0, \frac 1 3 \}$:
\begin{equation}
	\label{tab:ARS4_third_expl}
	\renewcommand{\arraystretch}{1.25}
	\text{ explicit: } \quad
	\begin{array}{c|cccc}
		0                             & 0                                                     & 0                                         & 0                       \\
		\frac{3 \gamma - 1}{6 \gamma} & \frac{3 \gamma - 1}{6 \gamma}                         & 0                                         & 0                       \\
		\frac{\gamma+1}{2}            & - \frac { 6\gamma^3 - 3\gamma^2 + 1} {2(3\gamma - 1)} & \frac{\gamma(3\gamma^2 + 1)}{3\gamma - 1} & 0                       \\
		\hline
		                              & 0                                                     & \frac{3 \gamma^2}{3\gamma^2 + 1}          & \frac{1}{3\gamma^2 + 1}
	\end{array}
	\quad \text{ implicit: } \quad
	\begin{array}{c|cccc}
		0                             & 0 & 0                                & 0                       \\
		\frac{3 \gamma - 1}{6 \gamma} & 0 & \frac{3 \gamma - 1}{6 \gamma}    & 0                       \\
		\frac{\gamma+1}{2}            & 0 & \gamma                           & \frac{1 - \gamma}{2}    \\
		\hline
		                              & 0 & \frac{3 \gamma^2}{3\gamma^2 + 1} & \frac{1}{3\gamma^2 + 1}
	\end{array}
\end{equation}

We now demonstrate that, using Butcher tableaux \eqref{tab:ARS4_third_expl},
a first order TVD scheme given by \eqref{eq:IMEXconvexStages} -- \eqref{eq:IMEXconvexUpdate} can be obtained.
First, we show the $L^\infty$ stability \eqref{def:LinftyStab}, i.e. $\|w^{n+1}\|_\infty \leq \|w^n\|_\infty$.
The idea of the proof lies in estimating the $L^\infty$ norm of each stage against $\|w^n\|_\infty$,
which will then be used to obtain the final estimate \eqref{def:LinftyStab}.
The proof is achieved by only applying the triangle and reverse triangle inequalities.
For clarity, we do not replace the coefficients of the tableaux
by their $\gamma$-dependent values yet.

The first stage in \eqref{eq:IMEXconvexStages} reduces to $w^{(1)} = w^n$ and we trivially obtain $\|w^{(1)}\|_\infty = \|w^n \|_\infty$.
Furthermore, $w^{(1)}$ is independent of $\theta_1$, and therefore we set $\theta_1 =1$,
as mentioned before.
Using definition \eqref{eq:IMEXconvexStages}, the second stage is given by
\begin{equation*}
	\label{eq:ARS4_second}
	w^{(2)}_j + a_{22}\impCFL \Delta_j^{(2)}= w^n_j - \lambda \left( (1- \theta_2) a_{22} \Delta_j^{n} + \theta_2 \tilde a_{21} \Delta_j^n\right).
\end{equation*}
Reformulating and collecting terms in $w_j^n, w_{j-1}^n$,
and using that $\tilde a_{21} = a_{22}$, we obtain the following expression
\begin{equation}
	\label{eq:ARS4_section_imp_as_LHS}
	w^{(2)}_j + a_{22}\impCFL \Delta_j^{(2)}= \left(1- \lambda a_{22} \right) w_j^n + \lambda a_{22} w_{j-1}^{n}.
\end{equation}
Note that the second stage, like the first one, is independent of $\theta_2$.
Thus, we set $\theta_2 = 1$,
which achieves the full contribution of the second stage of the IMEX-RK scheme.
For periodic boundary conditions and requiring $a_{22}$ and $1-\lambda a_{22}$ to be positive, we obtain the following estimates
\begin{align}
	\label{eq:LinfEstimate2}
	\begin{split}
		\|w^n\|_\infty & = \left( 1 -  \lambda a_{22}\right)\|w^n\|_\infty + \lambda a_{22} \|w^n\|_\infty  = \left( 1 -  \lambda a_{22}\right)~\max_j |w^n_j| + \lambda a_{22}~\max_j |w_{j-1}^n| \\
		& \geq \max_j \left|\left( 1 -  \lambda a_{22}\right)w_j^n + \lambda a_{22} w_{j-1}^n\right| = \max_j \left|w_j^n - \lambda a_{22} (w_{j}^n - w_{j-1}^n)\right|             \\
		& = \max_j \left|(1 + \impCFL a_{22})w_j^{(2)} - \impCFL a_{22} w_{j-1}^{(2)}\right|  \geq (1 + \impCFL a_{22})\|w^{(2)}\|_\infty - \impCFL a_{22}\|w^{(2)}\|_\infty        \\
		& = \|w^{(2)}\|_\infty.
	\end{split}
\end{align}
Further, from $a_{22}>0$ we immediately obtain a constraint on $\gamma$ given by  $\frac{3 \gamma -1}{6 \gamma} > 0$, i.e. $ \gamma > \frac{1}{3}$ or $\gamma < 0$.
Moreover, from $ 1 - \lambda a_{22} > 0$ we obtain a restriction on the time step given by $\lambda < \frac{6\gamma}{3\gamma -1}$.

Turning to the third stage of the scheme given in~\eqref{eq:IMEXconvexStages}, we obtain
\begin{equation}
	\label{eq:ARS4_third_before_simplification}
	\begin{split}
		w_j^{(3)} +(1-\theta_3) c_3 \impCFL\Delta_j^{(3)} = w_j^n &- \lambda \left((1-\theta_3)\tilde{c}_3 \Delta_j^n + \theta_3 \left(\tilde{a}_{31}\Delta_j^{(1)} + \tilde{a}_{32}\Delta_j^{(2)}\right)\right) \\
		&- \impCFL \theta_3 \left(a_{32}\Delta_j^{(2)} + a_{33}\Delta_j^{(3)}\right).
	\end{split}
\end{equation}
The next step consists in eliminating the terms in $\mu_\varepsilon w^{(2)}$ to avoid a $\varepsilon$ restriction on $\lambda$.
To that end, we invoke \eqref{eq:ARS4_section_imp_as_LHS} to
replace $-\impCFL \Delta_j^{(2)}$ in \eqref{eq:ARS4_third_before_simplification} by
\begin{equation}
	\label{eq:impcfl_second}
	-\impCFL \Delta_j^{(2)} = \frac{1}{a_{22}}(w_j^{(2)} - w_j^n) + \lambda \Delta_j^n,
\end{equation}
which yields, after some rearranging,
\begin{align}
	\label{eq:ARS4_third}
	\begin{split} 
		w^{(3)}_j + \impCFL \left((1-\theta_3) c_3 + \theta_3 a_{33}\right)\Delta_j^{(3)}
		=
		w_j^n
		&- \lambda\left((1-\theta_3)c_3 \Delta_j^n + \theta_3 \left(\tilde a_{31} \Delta_j^{n} + \tilde a_{32}\Delta ^{(2)}\right)\right) \\
		&+ \theta_3 a_{32}\left(\frac{1}{a_{22}}(w_j^{(2)} - w_j^n) + \lambda \Delta_j^{n}\right).
	\end{split}
\end{align}
Reformulating \eqref{eq:ARS4_third} to highlight the contribution of each step gives
\begin{align*}
	w^{(3)}_j + \impCFL \left(c_3 + \theta_3(a_{33} -c_3 )\right)\Delta_j^{(3)} =
	 & \left(1 - \lambda (1-\theta_3)c_3 - \theta_3 \lambda(\tilde a_{31}- a_{32}) - \frac{\theta_3 a_{32}}{a_{22}} \right)w_j^n         \\
	 & +\left(\lambda(1-\theta_3)c_3 + \theta_3 \lambda( \tilde a_{31}- a_{32})\right)w_{j-1}^n                                          \\
	 & + \theta_3 \left(\frac{ a_{32}}{a_{22}} - \lambda \tilde a_{32} \right)w_j^{(2)} + \theta_3~\lambda \tilde a_{32} w_{j-1}^{(2)}.
\end{align*}
To obtain the estimate $\|w^n\|_\infty \geq \|w^{(3)}\|_\infty$,
we proceed analogously to the proof for the second stage.
The proof follows the lines of \eqref{eq:LinfEstimate2} and is thus omitted.
Like in the second stage,
we find that the coefficients in front of all stages involving $w_j$ and $w_{j-1}$
have to be non-negative.
The non-negativity requirements on the coefficients of $w_{j}^{(2)}$ and $w_{j-1}^{(2)}$ result into the conditions
\begin{equation}
	\begin{aligned}
		\tilde a_{32}                                 & \geq 0 \iff \frac{\gamma(3\gamma^2 + 1)}{3\gamma - 1} \geq 0 \iff \gamma \leq 0 \text{ or } \gamma > \frac{1}{3}, \\
		\frac{a_{32}}{a_{22}} - \lambda \tilde a_{32} & \geq 0 \iff \lambda \leq \frac{6\gamma^2}{\gamma(3\gamma^2 + 1)} \text{ and } \gamma > \frac{1}{3}.
	\end{aligned}
\end{equation}
Further we find, thanks to the coefficient in front of $w_{j-1}^{n}$,
\begin{equation}
	\label{req:theta3}
	\begin{aligned}
		(1 - \theta_3) c_3 + \theta_3( \tilde{a_{31}} - a_{32} ) \geq 0 & \iff \theta_3 \frac{3 \gamma^2 ( \gamma + 1 )} { 3\gamma - 1 } \leq \frac{\gamma + 1} { 2 }\iff \theta_3 \leq \frac{3 \gamma - 1} { 6 \gamma^2 },
	\end{aligned}
\end{equation}
which yields a restriction on $\theta_3$ depending on $\gamma$.
Another estimate for $\theta_3$ can be obtained from the coefficient of $\Delta_{j}^{(3)}$ given by
\begin{equation*}
	\begin{aligned}
		c_3 + \theta_3(a_{33} -c_3 ) \geq 0 & \iff \gamma \theta_3 \leq \frac{\gamma + 1} { 2 }  \iff \theta_3 \leq \frac{\gamma + 1} { 2 \gamma }.
	\end{aligned}
\end{equation*}
We see that this condition on $\theta_3$ is less restrictive than
the one obtained from~\eqref{req:theta3} for all~$\gamma > \frac 1 3$.
Note that $\gamma < 0$ would yield a negative $\theta_3$,
which imposes $\gamma > 0$.
The largest value we can take for $\theta_3$ is therefore given by
\begin{equation*}
	\theta_3^\text{opt} = \frac{3 \gamma - 1} { 6 \gamma^2 },
\end{equation*}
and $\theta_3$ must satisfy $\theta_3 \leq \theta_3^\text{opt}$.
Note that the coefficient of $w_{j}^{n}$ is always non-negative.
The expression
\begin{equation}
	\label{req:lambda}
	1 - \lambda (1-\theta_3)c_3 - \theta_3 \lambda(\tilde a_{31}- a_{32}) - \frac{\theta_3 a_{32}}{c_2} \geq 0
\end{equation}
is always fulfilled if we set $\theta_3 = \theta_3^\text{opt}$.
Therefore, using $\theta_3^\text{opt}$ yields the maximal allowed input from the stages~\eqref{eq:IMEXstages} of the third order IMEX-RK scheme.
Otherwise,~\eqref{req:lambda} leads to another, more restrictive estimate for~$\lambda$.
Having obtained the estimate for stage three, we turn to the final update given by
\begin{align}
	\label{eq:UpdateTVD3}
	w_j^{n+1} + (1-\theta_4) \impCFL \Delta_j^{n+1}= w_j^n - \theta_4 (\lambda + \impCFL) b_2 \Delta_j^{(2)} - \theta_4 (\lambda + \impCFL) b_3 \Delta_j^{(3)} - (1-\theta_4) \lambda \Delta_j^n.
\end{align}
We again eliminate the terms in $\mu_\varepsilon \Delta_j^{(2)}$ and $\mu_\varepsilon \Delta_j^{(3)}$ in order to avoid a CFL restriction depending on $\varepsilon$.
In \eqref{eq:UpdateTVD3}, we replace $-\impCFL \Delta_j^{(2)}$ by \eqref{eq:impcfl_second} and $-\impCFL \Delta_j^{(3)}$, using \eqref{eq:ARS4_third}, by
\begin{align*}
	\begin{split}
		- \impCFL\Delta_j^{(3)} = \frac{1}{ c_3 + \theta_3 (a_{33}-c_3)} \left(w^{(3)}_j - w_j^n\right) &+ \frac{c_3}{c_3 + \theta_3 (a_{33}-c_3)}\lambda\left(\left(1 + \theta_3 (\tilde a_{31}-1) \right)\Delta_j^{n} + \theta_3 \tilde a_{32}\Delta_j ^{(2)}\right) \\
		&- \frac{\theta_3 c_3 a_{32}}{c_3 + \theta_3 (a_{33}-c_3)}\left(\frac{1}{a_{22}}(w_j^{(2)} - w_j^n) + \lambda \Delta_j^{n}\right).
	\end{split}
\end{align*}
Collecting the terms in front of the states $w_j, w_{j-1}$ at each stage and requiring the positivity of these coefficients gives the $L^\infty$ property $\|w^{n+1}\|_\infty \leq \|w^n\|_\infty$ following analogous steps to \eqref{eq:LinfEstimate2}.
The resulting requirements on the CFL condition $\lambda$ and the convex parameter $\theta_4$ are given in Lemma \ref{lem:Linf3}.

The TVD property can be obtained following the same line of computations as performed for the $L^\infty$ stability and is contained in the following result.
\begin{lemma}[$L^\infty$ stability, TVD property]
	\label{lem:Linf3}
	For periodic boundary conditions under the CFL condition
	\begin{equation*}
		\lambda \leq \frac{18 {\gamma^{3}} {\theta_4}-( 3 \gamma-1) ( 3 {\gamma^{2}}+1) }{( 3 \gamma-1) ( ( 6 {\gamma^{2}} + 1 ) \theta_4 - ( 3 {\gamma^{2}} + 1 ) ) },
	\end{equation*}
	the scheme consisting of the Butcher tableaux \eqref{tab:ARS4_third_expl} with the convex scheme given by the stages \eqref{eq:IMEXconvexUpdate} and the update \eqref{eq:IMEXconvexStages} is $L^\infty$ stable and TVD if the following conditions are fulfilled:
	\begin{equation}
		\gamma \geq \frac{\sqrt{3}}{3} , \quad \theta_1 = 1, \quad \theta_2 = 1, \quad \theta_3 = \frac{3 \gamma - 1} { 6 \gamma^2 }, \quad \theta_4 < \frac{( 3 \gamma-1) ( 3 {\gamma^{2}}+1) }{18 {\gamma^{3}}}.
	\end{equation}
\end{lemma}

\subsubsection{TVD conditions for arbitrary number of stages}
\label{sec:nth_order_convex}

In the spirit of the results from the third-order scheme,
we now seek a general framework on how to obtain TVD schemes with $s$ stages
using the IMEX formulation~\eqref{eq:IMEXconvexStages}--\eqref{eq:IMEXconvexUpdate}.
Since the proof follows steps analogous to the ones from Section~\ref{sec:third_order_convex},
we do not repeat the calculations and directly give the final result.
\begin{theorem}
	\label{theo:TVDconvex}
	Let $\tilde A, A \in \mathbb{R}^{s\times s}$, $\tilde b, b, \tilde c, c \in \mathbb{R}^s$
	define two Butcher tableaux~\eqref{tab:general} fulfilling~\eqref{cond:c} and the $p$-th order compatibility conditions.
	Let $\tilde b$ and $b$ coincide with the last rows of $\tilde A$ and $A$ respectively.
	For $k = 1, \ldots, s$ and $l = 1, \ldots, k-1$, we define
	\begin{equation*}
		\mathcal{A}_k = \theta_k a_{kk} + (1 - \theta_k)c_k, \quad
		\mathcal{\tilde A}_k = \theta_k a_{k1} + (1-\theta_k) \tilde c_k, \quad
		\mathcal{B}_{kl} = \frac{\theta_k a_{kl}}{\mathcal{A}_l}, \quad
		\mathcal{\tilde B}_{kl} = \theta_k \tilde a_{kl}.
	\end{equation*}
	In addition, we recursively define the following expressions:
	\begin{alignat*}{3}
		 & \mathcal{C}_k = \mathcal{\tilde A}_k - \sum_{l=2}^{k-1} \mathcal{B}_{kl} \mathcal{C}_l,             &                                                                                                                                & \qquad
		 &                                                                                                     & \mathcal{C}_{kl} = \mathcal{\tilde B}_{kl} - \sum_{r=l+1}^{k-1} \mathcal{B}_{kr} \mathcal{C}_{rl},                                       \\
		 & \mathcal{D}_k = 1 - \lambda \mathcal{\tilde A}_k - \sum_{l=2}^{k-1} \mathcal{B}_{kl} \mathcal{D}_l, &                                                                                                                                & \qquad
		 &                                                                                                     & \mathcal{D}_{kl} = \mathcal{ B}_{kl}- \lambda \mathcal{\tilde B}_{kl} - \sum_{r=l+1}^{k-1} \mathcal{B}_{kr} \mathcal{D}_{rl}.
	\end{alignat*}
	Then, with $\theta_1=1$ and under the following restrictions for $k = 2, \ldots, s$ and $l = 1, \ldots, k-1$,
	\begin{equation*}
		\mathcal{A}_k > 0,
		\quad \mathcal{C}_k \geq 0,
		\quad \mathcal{D}_k \geq 0,
		\quad \mathcal{C}_{kl} \geq 0,
		\quad \mathcal{D}_{kl} \geq 0.
	\end{equation*}
	the scheme consisting of the stages~\eqref{eq:IMEXconvexStages} and the update~\eqref{eq:IMEXconvexUpdate},
	combined with a TVD limiter, is $L^\infty$ stable and TVD under a CFL condition determined by $\lambda \geq 0$ where $\lambda$ does not depend on $\varepsilon$.
\end{theorem}

The result from Theorem \ref{theo:TVDconvex} concerns IMEX-RK schemes where the last rows of $\tilde{A}$ and $A$ respectively coincide with the weights $\tilde b$ and $b$.
This important class of schemes contain the so-called stiffly accurate (SA) IMEX-RK schemes, see e.g. \cite{ParRus2005} for details.
They are especially important in the context of low Mach applications, as discussed in the numerical results section.
However, the result of Theorem \ref{theo:TVDconvex} can be generalized to IMEX-RK schemes where the weights $\tilde b$ and $b$ do not coincide with the last rows of $\tilde{A}$ and $A$.
This generalization is briefly addressed in the following remark.

\begin{remark}
	For IMEX-RK schemes, where the weights $\tilde b$ and $b$ do not coincide with the respective last rows of $\tilde A$ and $A$ in the associated Butcher tableaux,
	we view the update~\eqref{eq:IMEXconvexUpdate} as an additional explicit $(s+1)$-th stage in \eqref{eq:IMEXconvexStages}.
	The matrices in the tableaux~\eqref{tab:general} are then of size $(s+1) \times (s+1)$ with the diagonal entry $a_{s+1,s+1} = 0$.
	The new weights $\tilde b$ and $b$ are then defined as the last rows of the new $\tilde A$ and $A$.
	Theorem~\ref{theo:TVDconvex} is then applied to yield the $L^\infty$ stability and the TVD property.
\end{remark}

We conclude this section with further remarks on an extension beyond CK type Butcher tableaux.

\begin{remark}
	We can prove the same result if the first column of $A$ allows for non-zero entries.
	The TVD conditions obtained when assuming that structure are given in~\ref{sec:appendix_non_CK_schemes}.
\end{remark}

\begin{remark}
	Applying this procedure to a two-stage, second-order scheme,
	we recover the ARS(2,2,2)-based schemes introduced in~\cite{DimLouMicVig2018}
	and studied more broadly in~\cite{MicTho2019}.
\end{remark}

\subsection{Fully discrete parachute scheme}
\label{sec:second_order_space}

To increase the resolution of the explicit spatial derivatives,
we provide a third-order reconstruction of the cell interface values $w_{j+1/2}$ such that the fully discrete scheme is $L^\infty$ stable and TVD.
We reconstruct the values $w_j^{(k)}$ using the neighbouring cell averages, see for instance~\cite{LeV1992}.
The reconstructed values $w_{j,-}^{(k)}$ and~$w_{j,+}^{(k)}$ at the inner interfaces of cell $C_j$ are then defined by
\begin{equation}
	\label{eq:space_reconstruction}
	w_{j,-}^{(k)} = w_j^{(k)} - \frac{\Delta x}{2} L\left( \sigma_{j + 1/2}^{(k)}, \sigma_{j - 1/2}^{(k)} \right), \quad
	w_{j,+}^{(k)} = w_j^{(k)} + \frac{\Delta x}{2} L\left( \sigma_{j - 1/2}^{(k)}, \sigma_{j + 1/2}^{(k)} \right),
\end{equation}
where $\sigma_{j + 1/2}^{(k)}$ denotes the slope between the values of $w_j^{(k)}$ and $w_{j+1}^{(k)}$ given by
\begin{equation*}
	\sigma_{j + 1/2}^{(k)} = \frac{ w_{j+1}^{(k)} - w_{j}^{(k)} } {\Delta x}.
\end{equation*}
The function $L( \sigma_L, \sigma_R )$ is a slope limiter
which ensures that the reconstructed values still satisfy $L^\infty$ property.
For a three-point stencil, the following estimate has to hold
\begin{equation}
	\label{def:TVD_Limiter}
	\min(|w_{j-1}^{(k)}|,|w_j^{(k)}|,|w_{j+1}^{(k)}|) \leq |w_{j,\pm}^{(k)}| \leq \max(|w_{j-1}^{(k)}|,|w_j^{(k)}|,|w_{j+1}^{(k)}|).
\end{equation}

Here, we use a TVD third-order space reconstruction satisfying~\eqref{def:TVD_Limiter} by following the limiting procedure introduced in~\cite{SchSeiTor2015}.
This procedure switches between the oscillatory unlimited third-order reconstruction and a third-order TVD limiter.
Switching to the TVD limiter is triggered by the appearance of a non-physical oscillation represented by a non-smooth extremum.

In the spirit of the reconstruction used to approximate the explicit derivatives, we could also increase the space accuracy of the implicit derivatives using TVD slope limiters.
Note that the slopes are determined in general by a non-linear function.
This implies an implicit approximation of the reconstructed values~\eqref{eq:space_reconstruction}.
Such computations usually lead to an iterative process
or a prediction-correction method,
and are therefore extremely costly.
We consider this increase in computational cost as too much in the sight of the actual gain in resolution.
Note that the IMEX-TVD scheme is overall only of first order.
Therefore, this is a loss of resolution in space we are willing to take to obtain a less costly TVD scheme.


\subsection{Choice of the free parameters}
\label{sec:numerics_parameters}

The TVD scheme of Section~\ref{sec:third_order_convex}, denoted by TVD$3$, contains some free parameters.
We now suggest optimal values of these free parameters.
To that end, we analyse the error produced by the schemes, as well as the CPU time taken, with respect to the free parameters.
This analysis will help us give some insights on how to optimally choose these parameters,
and on the trade-offs that must be made when making such choices.

We compare the IMEX$1$ scheme to the IMEX$3$, TVD$3$ and MOOD$3$ schemes.
Here, the IMEX$p$ scheme is the unlimited scheme of order $p$ in space and time.
Following this notation, the IMEX$1$ scheme is given by~\eqref{eq:firstOrder}
and the IMEX$3$ scheme corresponds to the Butcher tableaux~\eqref{tab:ARS4_third_expl}.
The MOOD$p$ scheme is obtained in the framework of Algorithm~\ref{algo:implicit_MOOD}, using the IMEX$p$ scheme as the high-order scheme and the TVD$p$ scheme as the parachute scheme.
For the linear multi-scale advection equation,
the detection criterion consists in taking the $L^\infty$ norm of the solution,
and thus we take $\Phi(w) = w$ in Algorithm~\ref{algo:implicit_MOOD}.
We also take $\xi = 0$ to eliminate all oscillations.

Since the scope of this section is to study the effect of the time discretisation on the precision and computational time of our schemes, we temporarily restrict ourselves to a first-order upwind space discretisation.
This ensures that only the effects of the time discretisation are monitored.

An exact smooth solution of the toy problem \eqref{eq:ScalarMultiscale} is used for the calibration of the free parameters $\gamma, \theta_3, \theta_4$ and $\lambda$ in the TVD3 scheme.
It is given by
\begin{equation}
	\label{eq:smooth_exact_solution}
	w^s(t,x) = 1 + \frac \varepsilon 2 \lp
	1 + \sin \left[ 2 \pi \varepsilon \lp x - \lp c_m + \frac {c_a} \varepsilon \rp t \rp \right]
	\rp,
\end{equation}
which represents a sine function of amplitude $\varepsilon$, transported with the velocity $c_m + \frac {c_a} \varepsilon$.
We set $\varepsilon = 0.1$ and we take $N = 400$ cells.
The conclusions of the forthcoming developments are unchanged if we consider other values of $\varepsilon$.
Indeed, taking a different~$\varepsilon$ would merely translate the curves without changing their relative positions.


\subsubsection{
	\texorpdfstring{Choice of $\gamma$ and $\theta_4$ in the TVD$3$ scheme}
	{Choice of gamma and theta4 in the TVD3 scheme}
}
\label{sec:choice_of_gamma_third_order}

\begin{figure}[!ht]

	\centering
	\makebox[0.95\linewidth][c]{

		\includegraphics{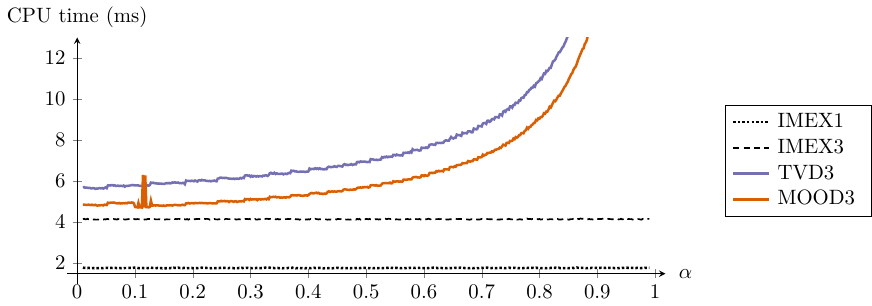}

	}

	\caption{
		CPU time (in milliseconds) with respect to the parameter $\alpha$,
		using $\gamma = \gamma^{\text{opt}} = \frac 2 3$,
		in the context of the test case presented in Section~\ref{sec:choice_of_gamma_third_order}.
	}

	\label{fig:CPU_time_wrt_alpha_third_order}

\end{figure}

For the TVD3 scheme, we have to set the values of $\theta_3$, $\theta_4$ and $\lambda$,
constrained by Lemma \ref{lem:Linf3}.
Ideally, we would like $\theta_3$, $\theta_4$ and $\lambda$ to be as large as possible.
By inspection, we note that the maximum value of $\theta_3$ is
$\theta_3^{\text{opt}} = \frac 3 8$, obtained for $\gamma^{\text{opt}} = \frac 2 3$.
The Butcher tableaux~\eqref{tab:ARS4_third_expl} then become
\begin{equation}
	\label{tab:ARS4_third_gamma_56}
	\renewcommand{\arraystretch}{1.25}
	\text{ explicit: }
	\begin{array}{c|cccc}
		0        & 0             & 0           & 0        \\
		\sfr 1 4 & \sfr 1 4      & 0           & 0        \\
		\sfr 5 6 & \sfr{-13}{18} & \sfr{14}{9} & 0        \\
		\hline
		         & 0             & \sfr 4 7    & \sfr 3 7
	\end{array},
	\text{\qquad \qquad \qquad}
	\text{ implicit: }
	\begin{array}{c|cccc}
		0        & 0 & 0        & 0        \\
		\sfr 1 4 & 0 & \sfr 1 4 & 0        \\
		\sfr 5 6 & 0 & \sfr 2 3 & \sfr 1 6 \\
		\hline
		         & 0 & \sfr 4 7 & \sfr 3 7
	\end{array}.
\end{equation}
Taking this value of $\gamma$ in Lemma \ref{lem:Linf3} yields the following bounds on $\theta_4$ and $\lambda$
\begin{equation}
	\label{eq:bounds_theta4_lambda}
	0 < \theta_4 < \frac 7 {16}
	\text{\quad and \quad}
	0 < \lambda < \frac { 7 - 16 \theta_4} { 7 - 11 \theta_4 }.
\end{equation}
We note that $\lambda$ is a decreasing function of $\theta_4$,
which implies that we are not able to use both a large~$\theta_4$ and a large $\lambda$.
There is a trade-off between the CFL condition $\lambda$, i.e. the CPU time, and the value of~$\theta$, i.e. the resolution of the scheme.
To quantify this balance between precision and CPU time, let us introduce $\alpha \in (0,1)$, to rewrite~\eqref{eq:bounds_theta4_lambda} as follows
\begin{equation}
	\label{eq:theta4_lambda_wrt_alpha}
	\theta_4 = \frac 7 {16} \alpha
	\text{\quad and \quad}
	\lambda = \frac { 1 - \alpha} { 1 - \frac{11}{16} \alpha }.
\end{equation}
We note that $\theta_4$ increases and~$\lambda$ decreases with increasing $\alpha$.
Making use of the formulation in terms of $\alpha$, we analyse the TVD3 scheme with
$\gamma = \gamma^{\text{opt}} = \frac 2 3$.
Note that a similar study was performed,
for the second-order case,
in~\cite{MicTho2019}.
We first display in Figure~\ref{fig:CPU_time_wrt_alpha_third_order} the CPU time with respect to $\alpha$ for the four schemes.
As expected, since the CFL condition becomes more restrictive, the CPU time increases with $\alpha$ for the TVD$3$ and the MOOD$3$ schemes.

\begin{figure}[!ht]

	\centering
	\makebox[0.95\linewidth][c]{

        \includegraphics{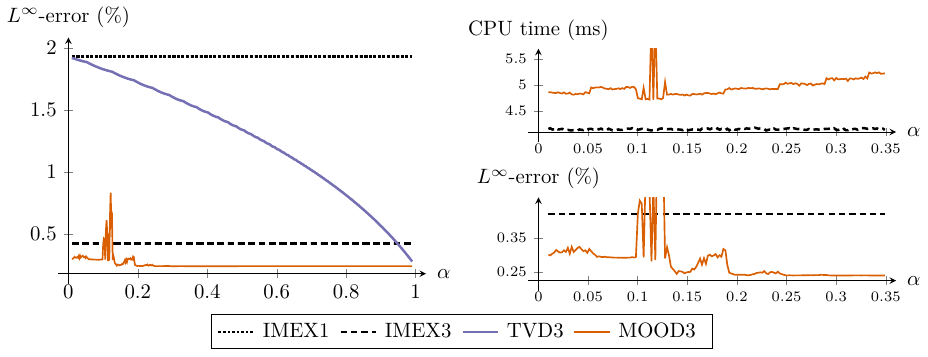}

	}

	\caption{
		$L^\infty$-error with respect to the parameter $\alpha$,
		using $\gamma = \gamma^{\text{opt}} = \frac 2 3$, $\theta_3 = \frac 3 8$
		and $\theta_4$, $\lambda$ given by relation \eqref{eq:theta4_lambda_wrt_alpha}.
		in the context of the test case presented in Section~\ref{sec:choice_of_gamma_third_order}.
		For $\alpha \in (0,0.35)$,
		the top right panel contains a zoom on the CPU time (data from Figure~\ref{fig:CPU_time_wrt_alpha_third_order})
		and the bottom right panel contains a zoom on the $L^\infty$-error (data from left panel).
	}

	\label{fig:Linf_error_wrt_alpha_third_order}

\end{figure}
Second, in the left panel of Figure~\ref{fig:Linf_error_wrt_alpha_third_order},
we display the $L^\infty$-error with respect to $\alpha$ for the four schemes under consideration.
As expected, we observe that it decreases with~$\alpha$ for the TVD$3$ scheme, since $\theta_4$ increases.
Third, in the right panel of Figure~\ref{fig:Linf_error_wrt_alpha_third_order},
we display a zoom on the CPU time and the $L^\infty$-error produced by the IMEX$3$ and MOOD$3$ schemes,
with respect to $0 < \alpha < 0.35$.
We observe that the error stabilizes around $\alpha = 0.3$, and that the CPU time increases monotonically with $\alpha$.
Therefore, taking $\alpha = \frac 1 3$ seems to be a good compromise between precision and computational time.
In the remainder of this article, we take
\begin{equation*}
	\label{eq:optimal_gamma_and_alpha}
	\gamma = \gamma^{\text{opt}} = \frac 2 3
	\text{\quad and \quad}
	\alpha = \alpha^{\text{opt}} = \frac 1 3,
\end{equation*}
which leads to the following values for $\theta_k, ~k=1,\dots,4$ and $\lambda$
associated to the Butcher tableaux \eqref{tab:ARS4_third_gamma_56}
\begin{equation*}
	\label{eq:optimal_thetas_lambda}
	\theta_1 = 1, \quad
	\theta_2 = 1, \quad
	\theta_3 = \theta_3^{\text{opt}} = \frac 3 8 = 0.375,
	\quad
	\theta_4 = \theta_4^{\text{opt}} = \frac 7 {48} \simeq 0.146
	\quad \text{and}\quad
	\lambda = \lambda^{\text{opt}} = \frac {32} {37} \simeq 0.865.
\end{equation*}


\subsubsection{
	\texorpdfstring{Numerical optimisation of larger Butcher tableaux}
	{Numerical optimisation of larger Butcher tableaux}
}
\label{sec:optimisaion_Butcher_tableaux}

To conclude this section, we mention a four-step third-order Butcher tableau yielding a TVD scheme.
To obtain this tableau, we have used the TVD inequalities from Theorem~\ref{theo:TVDconvex},
as well as the order conditions
as constraints in an optimisation problem.
Its objective is to maximize the value of $\lambda + \sum \theta$,
and its unknowns are the Butcher coefficients, the values of $\theta$, and $\lambda$.
We ran this optimization problem with many random initial conditions for the unknowns,
and we refined this random initialisation around values yielding a large value of the objective function.
In the end, we chose the solution where the value of the objective function was maximal,
under the additional constraint that~$\lambda \geq 0.5$, which is a standard CFL condition arising in fluid dynamics schemes.
We obtained $\lambda = 0.5471076190680170$, $\theta_1 = 1$, $\theta_2 = 1$, $\theta_3 = 1$, $\theta_4 = 0.5110907014643069$ and $\theta_5 = 0.4997722865197203$.
The Butcher tableau is given in~\ref{sec:appendix_TVD34}.
In the remainder of the paper, the scheme and its MOOD version will be referred to as TVD$3$($4$) and MOOD$3$($4$).


\section{Numerical results}
\label{sec:numerics}

In this section, we verify the properties of the numerical schemes that were developed in the previous sections.
The schemes are summarized in Table~\ref{tab:scheme_description}.
We first apply the TVD-MOOD strategy on the scalar linear problem from Section~\ref{sec:numerics_toy_problem},
before moving on to nonlinear systems of equations in Section~\ref{sec:isentropic_Euler}.
We consider the 2D isentropic Euler equations as an example of such a system.
We recall the two different types of CFL conditions used in the following.
The first one is given by an acoustic CFL condition
\begin{equation}
	\label{eq:def_nu_ac}
	\Delta t \leq \nu_{\text{ac}} \frac{\Delta x}{\underset{k}{\max}|\lambda_k|},
\end{equation}
which is restricted by the fastest wave speed $\max|\lambda_k|$.
The second one is a material CFL condition
\begin{equation}
	\label{eq:def_nu_mat}
	\Delta t \leq \nu_{\text{mat}} \frac{\Delta x}{|\lambda_u|},
\end{equation}
which is only constrained by the system velocity $|\lambda_u| \leq \max|\lambda_k|$ and thus allows larger time steps than the acoustic CFL condition.
For the toy problem, $\max|\lambda_k| = c_m + \frac{c_a}{\varepsilon}$
and $\lambda_u = c_m$ according to \eqref{eq:CFLac} and \eqref{eq:CFLmat}.

\begin{table}[!ht]
	\centering
	\begin{tabular}{clrlr}
		\toprule
		First-order IMEX \eqref{eq:firstOrder} & \multicolumn{2}{c}{Scheme from Section~\ref{sec:choice_of_gamma_third_order}} & \multicolumn{2}{c}{Scheme from Section~\ref{sec:optimisaion_Butcher_tableaux}}                                \\ \midrule
		                                       & unlimited:                                                                    & IMEX$3$                                                                        & unlimited:    & IMEX$3$($4$) \\
		parachute: IMEX$1$                     & parachute:                                                                    & TVD$3$                                                                         & parachute:    & TVD$3$($4$)  \\
		                                       & MOOD version:                                                                 & MOOD$3$                                                                        & MOOD version: & MOOD$3$($4$) \\ \bottomrule
	\end{tabular}
	\vspace{5mm}
	\caption{%
		Names of the schemes used in the numerical simulations.
	}
	\label{tab:scheme_description}
\end{table}


\subsection{Linear multi-scale advection equation}
\label{sec:numerics_toy_problem}
Since oscillations typically appear around discontinuities, we consider the following exact solution of the initial value problem~\eqref{eq:ScalarMultiscale}
\begin{equation}
	\label{eq:discontinuous_exact_solution}
	w(t,x) =
	\begin{dcases}
		1 + \varepsilon & \text{ if } \frac 1 4 <
		\left(
		\frac{\lp x - \lp c_m + \frac {c_a} \varepsilon \rp t \rp}{c_m + \frac{c_a}{\varepsilon}}  -
		\left\lfloor \frac{\lp x - \lp c_m + \frac {c_a} \varepsilon \rp t \rp}{c_m + \frac{c_a}{\varepsilon}} \right\rfloor
		\right)
		< \frac 3 4,                               \\
		1               & \text{ otherwise,}
	\end{dcases}
\end{equation}
which represents a rectangular bump located at $\left(\frac 1 4 ( c_m + \frac {c_a} \varepsilon ), \frac 3 4 ( c_m + \frac {c_a} \varepsilon )\right)$ of amplitude $\varepsilon$.
It is transported with the velocity $c_m + \frac {c_a} \varepsilon$
and the initial condition is given by the solution at time $t=0$.
The computational domain is given by $(0, c_m + \frac{c_a}{\varepsilon})$ with final time $T_f = 1$, which corresponds to one full cycle with periodic boundary conditions.
We take $c_m = 1$ and $c_a = 1$.

Since we consider a linear advection problem that respects a maximum principle, the detection criterion~$\Phi$ in the MOOD Algorithm \ref{algo:implicit_MOOD} is given by the $L^\infty$ norm of the solution
\begin{equation}
	\label{eq:detection_criterion_toy_problem}
	\Phi(w) = w.
\end{equation}
The hierarchy of the MOOD procedure consists in a high-order scheme given by
the third-order IMEX3 or IMEX3(4) schemes from Section~\ref{sec:choice_of_gamma_third_order} and
Section~\ref{sec:optimisaion_Butcher_tableaux}
with centred differences in the implicitly treated space derivative,
whereas the parachute scheme is given by the first-order IMEX scheme TVD3 or TVD3(4) with upwind discretization for all derivatives.

We compare our results against the $L$-stable ARS(2,3,3) scheme, reported in~\cite{AscRuuSpi1997}
or~\cite{ParRus2001} and recalled in~\ref{sec:appendix_ARS223}.
Note that the ARS(2,3,3) scheme can be written in the form of \eqref{tab:ARS4_third_expl}
with $\gamma = \frac{3 - \sqrt{3}}{6}$, and therefore it falls within the framework of CK schemes that were used in Section \ref{sec:third_order_convex}.
However, this value of~$\gamma$ does not satisfy the requirement of Lemma \ref{lem:Linf3}.
Consequently, within our framework, we cannot prove the existence of a convex combination, starting with the ARS(2,3,3) scheme, that leads to a first order TVD and~$L^\infty$ stable scheme with a material CFL condition \eqref{eq:def_nu_mat}.
In the following, we display the numerical results for the IMEX$1$, ARS(2,3,3), IMEX$3$, TVD$3$, MOOD$3$ and MOOD$3$($4$) schemes.


\subsubsection{Maximum principle}
\label{sec:discontinuous_approximation}

\newcommand{\Lover}{{L^1_{\text{o}}}}

\begin{figure}[!b]

	\centering
	\makebox[0.95\textwidth][c]{

		\includegraphics{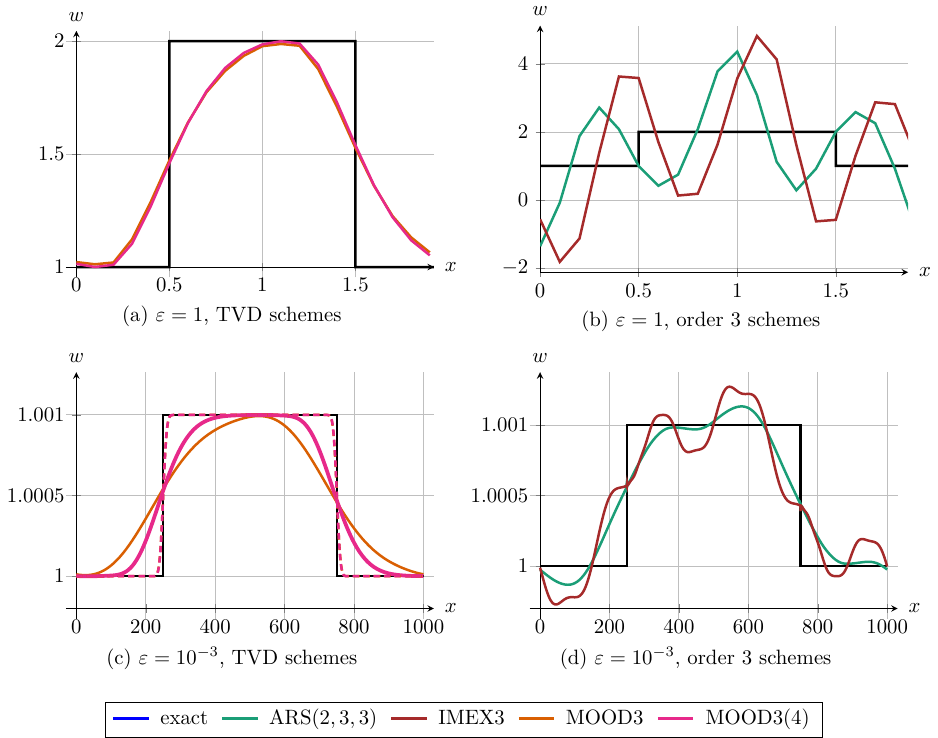}

	}

	\caption{
		Discontinuous solution~\eqref{eq:discontinuous_exact_solution} of the linear advection problem at time $T_f = 1$
		with $\Delta x = 0.1$ for $\varepsilon = 1$ (top) and $\varepsilon = 10^{-3}$ (bottom).
        The dashed lines denote the acoustic CFL condition \eqref{eq:def_nu_ac} with $\nu_{\text{ac}} = 0.5$,
        and the solid lines use the material CFL condition \eqref{eq:def_nu_mat} with $\nu_{\text{mat}} = 0.5$.
	}

	\label{fig:step_order_3_solutions}

\end{figure}

We study the behaviour of the above-mentioned schemes with respect to the maximum principle for the discontinuous solution~\eqref{eq:discontinuous_exact_solution}.
First, in Figure~\ref{fig:step_order_3_solutions}, the numerical approximations of the discontinuous solution \eqref{eq:discontinuous_exact_solution}
for $\varepsilon = 1$ and $\varepsilon = 10^{-3}$ are displayed.
For $\varepsilon=1$, the solutions are computed with a material CFL condition given by $\nu_{\text{mat}} = 0.5$.
For $\varepsilon=10^{-3}$, the results are given in addition using an acoustic CFL condition given by $\nu_{\text{ac}} = 0.5$.
These CFL conditions respectively correspond to time steps of $\Delta t_{\text{ac}} \simeq 10^{-5}$ and $\Delta t_{\text{mat}} = 10^{-2}$.
In both cases, we take $\Delta x = 0.1$, which corresponds to~$20$ cells for $\varepsilon=1$ and $10010$ cells for $\varepsilon=10^{-3}$.

Since the solution contains only fast travelling waves, using a material time step leads to diffused wave fronts.
However, for the acoustic time step, they are captured accurately since the CFL condition is limited by these fast waves.
We note that the purely third order schemes IMEX3 and ARS($2,3,3$) are oscillatory and violate the maximum principle even for $\varepsilon = 1$, even though the ARS($2,3,3$) scheme is $L$-stable.
For $\varepsilon = 10^{-3}$, the unlimited third-order schemes are not in-bounds when using the material CFL condition.
Their oscillations are not too large thanks to the diffusion from the upwind treatment of the implicit derivative.
However, with the acoustic CFL condition, the centred differences on the implicit derivative produce extremely large oscillations which are not displayed in Figure~\ref{fig:step_order_3_solutions}.
Conversely, the MOOD schemes are in-bounds for both choices of CFL condition and both $\varepsilon$ regimes.
Note that the IMEX$3$($4$) scheme is more stable and precise than the IMEX$3$ scheme,
which results in a better precision for the MOOD$3$($4$) scheme compared to the MOOD$3$ scheme.
Indeed, the lack of stability of the IMEX$3$ scheme, especially for small $\varepsilon$,
triggers the parachute scheme more often to ensure the maximum principle.

Next, we compute the error on the numerical solution.
A standard error computation in the context of finite volume schemes consists in using the $L^1$ norm, defined by
\begin{equation*}
	\label{eq:L1_norm}
	\| w^n \|_1 = \frac 1 {\Delta x} \sum_j |w_j^n|.
\end{equation*}
However, the above norm only measures the average deviation between the exact solution and the numerical approximation.
Since we seek to measure of the violation of the maximum principle, we need to take into account the impact of overshoots and undershoots on the error.
Therefore, we propose to use the following quasinorm
\begin{equation*}
	\label{eq:L1o_norm}
	\| w^n \|_{\Lover} = \frac 1 {\Delta x} \sum_j \lp |w_j^n| + \max_{m \leq n} \left[ \lp \max_j w_j^m - \min_j w_j^m \rp -
		\lp \max_j w_j^0 - \min_j w_j^0 \rp \right] \rp,
\end{equation*}
which adds the error of the over- and undershoots to the $L^1$ error.

\begin{figure}[!ht]

	\centering
	\makebox[0.95\textwidth][c]{

		\includegraphics{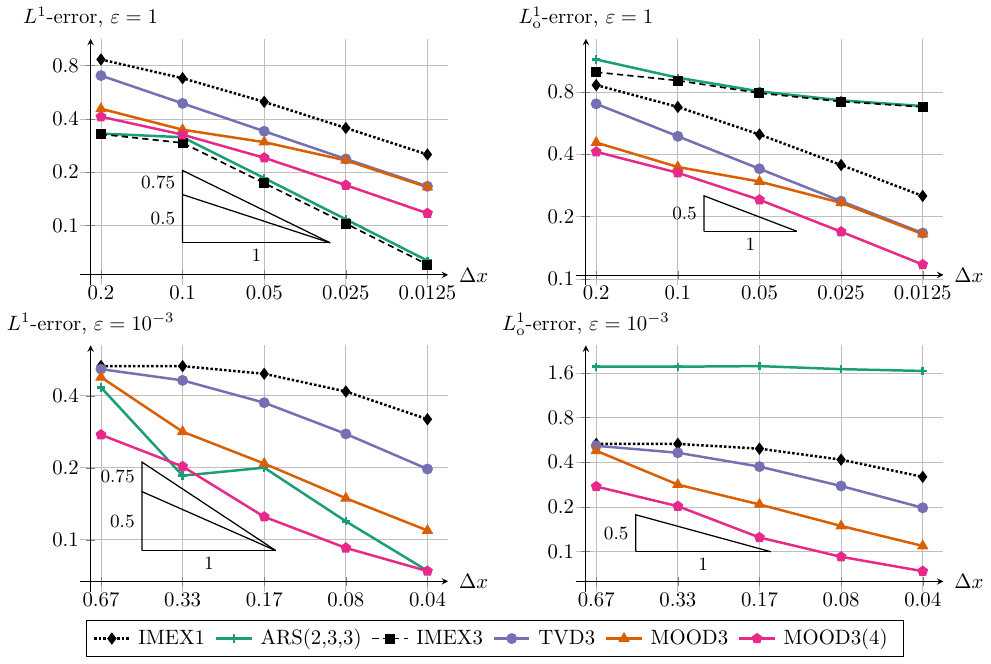}

	}

	\caption{
		Error lines in $L^1$ norm (left) and $\Lover$ quasinorm (right)
		for the discontinuous solution~\eqref{eq:discontinuous_exact_solution} with $\varepsilon = 1$ (top) and $\varepsilon = 10^{-3}$ (bottom).
	}

	\label{fig:step_order_3_error_lines}

\end{figure}

In Figure~\ref{fig:step_order_3_error_lines}, we report the error in the $L^1$ norm and in the $\Lover$ quasinorm produced by the above-mentioned schemes for $\varepsilon = 1$ and $\varepsilon = 10^{-3}$
with $\nu_{\text{mat}} = 0.5$,
corresponding to $\Delta t = 0.01$.
First, we observe that the theoretical order of convergence for a higher order scheme is not reached since we obtain an accuracy up to order $\frac 1 2$ for the IMEX$1$, TVD$3$, MOOD$3$ and MOOD$3$($4$) schemes,
and up to order $\frac 3 4$ for the ARS(2,3,3) and IMEX$3$ schemes.
This is due to the fact that we approximate a discontinuous solution
where the numerical diffusion of the schemes considerably reduces the order of convergence, see for instance~\cite{LeV1992}.
Second, when taking the over- and undershoots into account, i.e. considering the error in the $\Lover$ quasinorm,
we observe that the IMEX$1$, TVD$3$, MOOD$3$ and MOOD$3$($4$) schemes show comparable experimental order of convergence (EOC) as in the $L^1$ norm.
This was to be expected since they are constructed such that over- or undershoots are avoided.
However, for the ARS(2,3,3) scheme, the error in the $\Lover$ norm is roughly constant even when an increasing number of cells is considered.
This means that the improvement of the $L^1$ error is almost exactly compensated by an increase of the over- and undershoots in absolute magnitude.
Therefore, even taking an increasing number of cells is not enough
to stabilise the ARS(2,3,3) scheme.

\subsubsection{Advantages of the TVD-MOOD approach}
\label{sec:perks}

We continue our discussion by addressing the flexibility of our schemes with respect to the choice of time stepping compared to $L$-stable and SSP IMEX schemes from the literature.
Further, we show that our first order IMEX-TVD schemes, overall, improve the space-time errors compared to the classical IMEX$1$ first-order backward forward Euler scheme.


\paragraph{Flexibility on the time step}
In the following, we compare the MOOD$3$($4$) scheme with higher order schemes from the literature.
Namely, we consider the aforementioned ARS(2,3,3) scheme from~\cite{AscRuuSpi1997},
as well as a more recent SSP-IMEX scheme from~\cite{ConGotGraSha2017} which we refer to as CGGS$3$.
In contrast to our MOOD$3$($4$) scheme, the ARS(2,3,3) and CGGS$3$ schemes are restricted by an acoustic CFL condition to be $L^\infty$ stable.
This results in a $\varepsilon$-dependent CFL restriction $\lambda \leq 0.9\varepsilon$, which corresponds to $\nu_{ac} = 0.9$ and $\nu_{mat} = 0.009$.
Note that our MOOD$3$($4$) scheme allows $\lambda < \lambda^{opt} \simeq 0.547$ corresponding to $\nu_{ac} = 54.7$ and $\nu_{mat} = 0.547$.
In Table \ref{tab:flexibility}, the CPU time and the $L^1$ error
for the discontinuous solution~\eqref{eq:discontinuous_exact_solution} with $N = 4000$ and $\varepsilon = 10^{-3}$ for the MOOD$3$($4$), CGGS$3$ and ARS(2,3,3) schemes are displayed.
Recall that the time step is $\Delta t = \Delta x \frac{\lambda}{1 + 1/\varepsilon} \simeq 4 \lambda$.
Applying an acoustic CFL condition gives comparable errors and CPU times for all schemes,
leading to a resolution of all waves.
However, in the case of our MOOD$3$($4$) scheme, the scheme is still stable for larger values of $\lambda$.
This can be especially advantageous if the resolution of certain waves can be neglected, which reduces computational time.
Examples of this behaviour include the results with a material time step given in Figure~\ref{fig:step_order_3_solutions},
or the approximation of material waves in the context of the isentropic Euler equations in Figure~\ref{fig:Euler_material_wave}.
Also, note that using high-order schemes such as~ARS(2,3,3) or CGGS$3$ enforces a time step that vanishes if~$\varepsilon$ tends to zero,
leading to long computational times especially when flow phenomena are monitored over long time intervals.

\begin{table}[!ht]
	\centering
	\renewcommand{\arraystretch}{1.1}
	\begin{tabular}{@{}lllllll@{}} \toprule
		             &  & $\lambda$                                     &  & CPU time (s) &  & $L^1$ error \\ \midrule
		MOOD$3$($4$) &  & $\lambda = \lambda^{\text{opt}} \simeq 0.548$ &  & 0.0101       &  & 0.217       \\
		             &  & $\lambda = 250 \varepsilon = 0.25$            &  & 0.0222       &  & 0.111       \\
		             &  & $\lambda = 50 \varepsilon = 0.05$             &  & 0.0953       &  & 0.0591      \\
		             &  & $\lambda = 10 \varepsilon = 0.01$             &  & 0.659        &  & 0.0488      \\
		             &  & $\lambda = 2 \varepsilon = 0.0002$            &  & 1.63         &  & 0.0253      \\
		             &  & $\lambda = 0.9 \varepsilon = 0.0009$          &  & 3.64         &  & 0.0253      \\\midrule
		CGGS$3$      &  & $\lambda = 0.9 \varepsilon = 0.0009$          &  & 1.25         &  & 0.0253      \\\midrule
		ARS(2,3,3)   &  & $\lambda = 0.9 \varepsilon = 0.0009$          &  & 1.17         &  & 0.0253      \\ \bottomrule
	\end{tabular}
	\vspace{\smallskipamount}
	\caption{%
	CPU times and $L^1$ errors for discontinuous solution~\eqref{eq:discontinuous_exact_solution} with $N = 4000$ discretisation points and $\varepsilon=10^{-3}$,
	using the  MOOD$3$($4$), ARS(2,3,3) and CGGS$3$.%
	}
	\label{tab:flexibility}
\end{table}


\paragraph{Space-time error of the TVD-MOOD approach}
In this paragraph, we study the impact of using the TVD$3$($4$) scheme as a parachute scheme instead of the usual IMEX$1$ backward forward Euler scheme.
For fairness, we compare the TVD$3$($4$) scheme to the four-stage IMEX$1$($4$) version of the IMEX$1$ scheme.
In addition, since we wish to compare the time discretisations, we only consider first-order upwind space discretisations in this paragraph.
Note that both the TVD$3$($4$) and IMEX$1$($4$) schemes are first-order accurate.
An important point to quantify is the reduction of diffusion reflected in the error made by using our TVD schemes as parachute schemes in the MOOD Algorithm \ref{algo:implicit_MOOD}.
Therefore, we introduce space-time errors, which measure the average and maximum errors between the numerical and exact solution caused by diffusion over the whole computational time.
Note that, in contrast, $L^1$ errors are usually computed at the final time only.
The space-time errors are defined as
\begin{equation*}
	\label{eq:space_time_error}
	\begin{aligned}
		e_{\text{ST}}^{\text{mean}} & = \frac 1 {n_t} \sum_{n = 1}^{n_t} \left[ \left( \max_j (w_{\text{ex}})_j^n - \min_j (w_{\text{ex}})_j^n \right) - \left( \max_j w_j^n - \min_j w_j^n \right) \right], \\
		e_{\text{ST}}^{\text{max}}  & = \max_{1 \leq n \leq n_t} \left[ \left( \max_j (w_{\text{ex}})_j^n - \min_j (w_{\text{ex}})_j^n \right) - \left( \max_j w_j^n - \min_j w_j^n \right) \right],
	\end{aligned}
\end{equation*}
where $n_t$ is the number of time iterations.
In Table~\ref{tab:space_time_errors}, we report the values of $e_{\text{ST}}^{\text{mean}}$ and $e_{\text{ST}}^{\text{max}}$ for the discontinuous solution \eqref{eq:discontinuous_exact_solution} with respect to $\varepsilon \in \{1, 10^{-1}, 10^{-2}, 10^{-3}\}$.
We take $\Delta x = 0.1$, leading to $N \simeq 20 / \varepsilon$ cells since the size of the computational domain depends on $\varepsilon$.
We set $\nu_{mat} = 0.5$, which corresponds to $\nu_{ac} = \nu_{mat} (c_m + c_a / \varepsilon) / c_m = 0.5 + 0.5 / \varepsilon$ and $\Delta t = 0.01$.
We note that, for the mean and maximum space-time errors, using the TVD$3$($4$) scheme as a parachute scheme lowers the error by a factor of up to $100$ for small values of $\varepsilon$, compared to the IMEX$1$($4$) scheme.
This is due to the fact that the TVD$3$($4$) scheme is much less diffusive than the IMEX$1$($4$) scheme
for small values of $\varepsilon$, and thus the approximate solution stays closer to the exact solution when the parachute scheme is triggered,
rather than being diffused away.

\begin{table}[!ht]
	\makebox[\textwidth][c]{
		\begin{tabular}{@{}ccccccccc@{}} \toprule
			$\varepsilon$ &  & \multicolumn{3}{c}{mean spacetime error} &                              & \multicolumn{3}{c}{maximum spacetime error}                                                                            \\ \cmidrule{1-9}
			              &  & \faParachuteBox: IMEX$1$($4$)            & \faParachuteBox: TVD$3$($4$) & ratio                                       &  & \faParachuteBox: IMEX$1$($4$) & \faParachuteBox: TVD$3$($4$) & ratio  \\ \cmidrule(lr){3-5} \cmidrule(lr){7-9}
			$1.0$         &  & $2.29 \times 10^{-1}$                    & $2.29 \times 10^{-1}$        & $1.00$                                      &  & $5.20 \times 10^{-1}$         & $5.20 \times 10^{-1}$        & $1.00$ \\
			$10^{-1}$     &  & $3.41 \times 10^{-3}$                    & $3.91 \times 10^{-3}$        & $0.87$                                      &  & $1.80 \times 10^{-4}$         & $1.99 \times 10^{-5}$        & $0.90$ \\
			$10^{-2}$     &  & $3.03 \times 10^{-7}$                    & $1.46 \times 10^{-8}$        & $20.8$                                      &  & $1.41 \times 10^{-6}$         & $1.30 \times 10^{-7}$        & $10.8$ \\
			$10^{-3}$     &  & $2.93 \times 10^{-6}$                    & $3.21 \times 10^{-8}$        & $91.3$                                      &  & $2.37 \times 10^{-5}$         & $2.46 \times 10^{-7}$        & $96.3$ \\ \bottomrule
		\end{tabular}
	}
	\vspace{\smallskipamount}
	\caption{%
	Space-time error with respect to $\varepsilon$ using the TVD3(4) and IMEX1(4) scheme as a parachute for the discontinuous solution \eqref{eq:discontinuous_exact_solution}.
	}
	\label{tab:space_time_errors}
\end{table}


\subsection{Isentropic Euler equations}
\label{sec:isentropic_Euler}

In this section, we apply the TVD-MOOD strategy given by Algorithm \ref{algo:implicit_MOOD} to the isentropic Euler equations.
They are governed in a non-dimensional formulation by
\begin{equation}
	\label{eq:isentropic_euler}
	\begin{dcases}
		\partial_t \rho + \nabla \cdot (\rho \bm u) = 0, \\
		\partial_t (\rho \bm u) + \nabla \cdot (\rho \bm u \otimes \bm u) + \frac 1 {M^2} \nabla p(\rho) = 0,
	\end{dcases}
\end{equation}
where $\rho(x,t) > 0$ denotes the density and $\bm u(x,t)$ the velocity field.
Assuming an ideal gas, the pressure law is given by~$p(\rho) = \rho^\gamma$, where $\gamma \geq 1$ is the ratio of specific heats.
Due to the non-dimensional formulation, the pressure gradient is scaled by the Mach number squared,
which is given by the ratio between fluid velocity~$\|\bm u\|$ and sound speed $c$.
This leads to Mach number-dependent acoustic wave speeds $\lambda^\pm = u \pm \frac{c}{M}$ with $c^2 = \partial_\rho p$.
For small Mach numbers, they tend to infinity, and are therefore integrated implicitly.
This results in the following splitting of the flux
\begin{equation}
	\partial_t w + \nabla \cdot f_e(w) + \nabla \cdot f_i(w) = 0,
\end{equation}
with the state vector $w = (\rho, \rho \bm u)^T$, the explicitly treated flux $f_e(w) = (0, \rho \bm u \otimes \bm u)^T$
and the implicitly treated flux $f_i(w) = (\rho \bm u, \frac{1}{M^2} p I)^T$.
The reader is referred, for instance, to~\cite{KlaMaj1981,Del2010}
for a description of the low Mach number limit,
and to \cite{DegTan2011,BosRusSca2018,ThoZenPupKli2020} for more information on the so-called asymptotic-preserving property,
which ensures the consistency of the numerical scheme with the incompressible Euler equations in the singular Mach number limit.

Regarding the space discretisation, we consider two schemes given by
a RS-IMEX scheme in the spirit of \cite{ZeiSchKaiBecLukNoe2020,LukPupTho2022}
and the TVD IMEX scheme from \cite{DimLouVig2017}.
Both of them have a Mach number independent CFL condition.
The RS-IMEX scheme is obtained by linearising the pressure against a constant reference density, as it appears in the low Mach limit,
and the implicitly treated flux is discretised with centred differences.
This prevents the scheme from being TVD, but is consistent with the low Mach number limit.
		An extension to non-stationary reference solutions can be obtained following a prediction-correction approach given in~\cite{ParMun2005}.
In the scheme from \cite{DimLouVig2017}, all fluxes are discretized in an upwind fashion, which leads to a TVD scheme, but with a numerical viscosity depending on the Mach number, leading to possibly non-feasible numerical solutions in the low Mach number limit.
We thus propose the following MOOD hierarchy.
First, the high order scheme is the third order IMEX$3$($4$) scheme with the RS-IMEX discretisation.
Then, we introduce an intermediate stage with a TVD$3$($4$) integration and RS-IMEX discretization.
It has a Mach number-independent diffusion but is overall not TVD due to the centred differences.
Finally, the parachute scheme is the TVD$3$($4$) scheme with the upwind discretisation from~\cite{DimLouVig2017},
which is TVD but very diffusive for low Mach numbers, and should only be triggered rarely.

For the definition of the MOOD criterion \eqref{eq:detection_criterion_toy_problem}, we follow~\cite{DimLouMicVig2018} and use the Riemann invariants to detect oscillations.
This choice is motivated by the fact that, according to~\cite{SmoJoh1969}, at least one of the Riemann invariants satisfies a maximum principle in a 1D Riemann problem.
It can be extended to the~2D Cartesian framework by using the normal velocity to define the Riemann invariants at each edge of the mesh.
The criterion is then given by
\begin{equation}
	\label{eq:MOODcrit}
	\Phi_\pm(W) = \bm u \cdot \bm n \mp \frac{1} {M} \frac 2 {\gamma - 1} \sqrt{\gamma \rho^{\gamma - 1}},
\end{equation}
where $\boldsymbol n$ is the outward-pointing normal vector.
The detection criterion
in Algorithm~\ref{algo:implicit_MOOD} is then time-dependent and given by $\Phi = \max(\Phi_-, \Phi_+)$.
Recall from Algorithm~\ref{algo:implicit_MOOD} that $\xi \in [0, 1]$
controls the permitted oscillations in the MOOD solution.
Choosing a small $\xi$ allows the presence of small oscillations in the MOOD solution.
However, setting $\xi = 0$ would strictly enforce the detection criterion, and the parachute scheme would be triggered too often, leading to a very diffused solution especially in low Mach number regimes.
Therefore, we set $\xi = \frac 1 {100}$ moderately small for all the following numerical experiments.
Due to the splitting, the stability of the numerical scheme solely depends on the fluid velocity.
The CFL condition in normal direction is thus given by
\begin{equation}
	\label{eq:time_step_Euler}
	\Delta t \leq \nu_{\text{mat}} \frac {\Delta x} {2 \max|\bm u \cdot \bm n|},
\end{equation}
which does not depend on the Mach number $M$.
The acoustic CFL condition \eqref{eq:def_nu_ac} is restricted by $\max|\lambda^\pm|$ and enforces a vanishing time step when $M$ tends to 0.
In the following numerical tests, we focus on the two-dimensional setting on Cartesian grids.

\subsubsection{Riemann problems}
\label{sec:2D_RP}

We first apply the scheme on three Riemann problems.
The first one only concerns the propagation of acoustic waves,
the second one contains a slow moving material wave, while the third one concerns the correct capturing of shock waves in a circular setting.
For each problem, if not stated otherwise, we provide a reference solution computed using the IMEX$1$ scheme on a fine grid.
Since the jump in the initial condition is only in $x$- or in radial direction,
these tests mimic one-dimensional Riemann problems.
Therefore, we consider for the first two Riemann problems, a mesh with $100$ cells in the $x$-direction and $3$ cells in the $y$-direction, $\gamma = 1.4$ and the computational domain is given by $[0, 2] \times [0, 1]$, equipped with Neumann boundary conditions.

\begin{figure}[!ht]

	\centering
	\makebox[0.95\linewidth][c]{

		\includegraphics{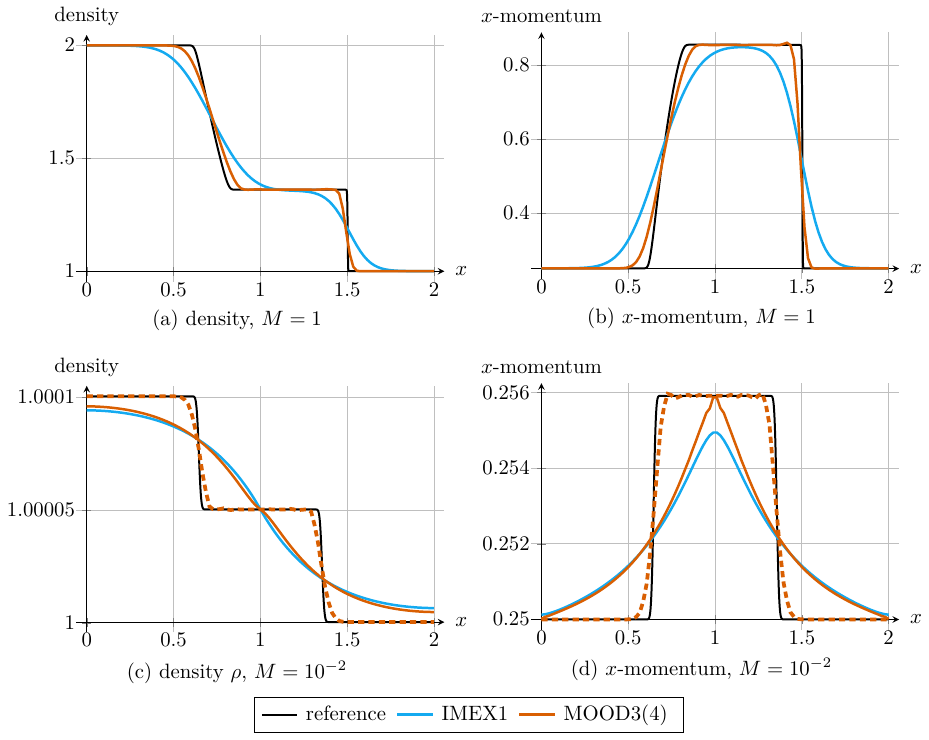}

	}

	\caption{Acoustic waves:
		Approximation of the solution to the Riemann problem~\eqref{eq:Euler_Riemann_problem}
		at time $T_f = 0.3M$ on a $100 \times 3$ mesh,
		with $M = 1$ (top)
		and $M = 10^{-2}$ (bottom).
		For $M = 10^{-2}$, we report the results with material $\nu_{\text{mat}} = 0.5$ (solid lines) and acoustic $\nu_{\text{ac}} = 0.5$ (dashed lines) CFL restrictions.
	}

	\label{fig:Euler_RP}

\end{figure}

\paragraph{Acoustic waves}
The initial data of the Riemann problem is given by
\begin{equation}
	\label{eq:Euler_Riemann_problem}
	\rho(x,y,0) =
	\begin{dcases}
		1 + M^2 & \text{ if } x < 1, \\
		1       & \text{ otherwise,}
	\end{dcases}, \quad \bm u(x,y,0) = 0.
\end{equation}
The approximations are depicted at the final time $T_f = 0.3 M$ in Figure~\ref{fig:Euler_RP} for $M =1$ and $M = 10^{-2}$,
with $\nu_{\text{ac}} = 0.5$ and $\nu_{\text{mat}} = 0.5$.
For $M = 10^{-2}$, these CFL conditions respectively correspond to
$\Delta t_{\text{ac}} \simeq 8.43 \times 10^{-5}$ and $\Delta t_{\text{mat}} = 2 \times 10^{-2} > 0.3 M$ for $M = 10^{-2}$.

We note that, using the acoustic CFL condition, the shock and rarefaction waves are captured accurately with the correct propagation speed for both Mach numbers.
Applying a material CFL condition leads to a stable scheme for $M = 10^{-2}$.
As expected for larger time steps, the acoustic waves are diffused but the MOOD3(4) scheme still gives the correct maximal momentum.
Utilizing such large time steps allows to get a correct description of the solution shape using only one time step, compared to $36$ time steps with an acoustic CFL condition.

Furthermore, for $M=1$ and $M=10^{-2}$,
the MOOD procedure is activated respectively $51\%$ and $30\%$ of the time iterations using acoustic time stepping.
This observation is explained by the fact that the IMEX$3$($4$) scheme is quite oscillatory,
which have to be corrected by the MOOD procedure.
The share of iterations where the MOOD criterion \eqref{eq:MOODcrit} was violated
could be lowered by basing the MOOD procedure on a less oscillatory third-order IMEX scheme than the IMEX3(4) scheme, as will be done in Section~\ref{sec:Euler_vortex}.
Note that, as expected, setting a non zero value of~$\xi$ in the MOOD procedure gives rise to marginal oscillations which do not interfere with the quality of the solution.

\paragraph{Material wave}
\label{sec:2D_shear}

\begin{figure}[!ht]

	\centering
	\makebox[0.95\linewidth][c]{

		\includegraphics{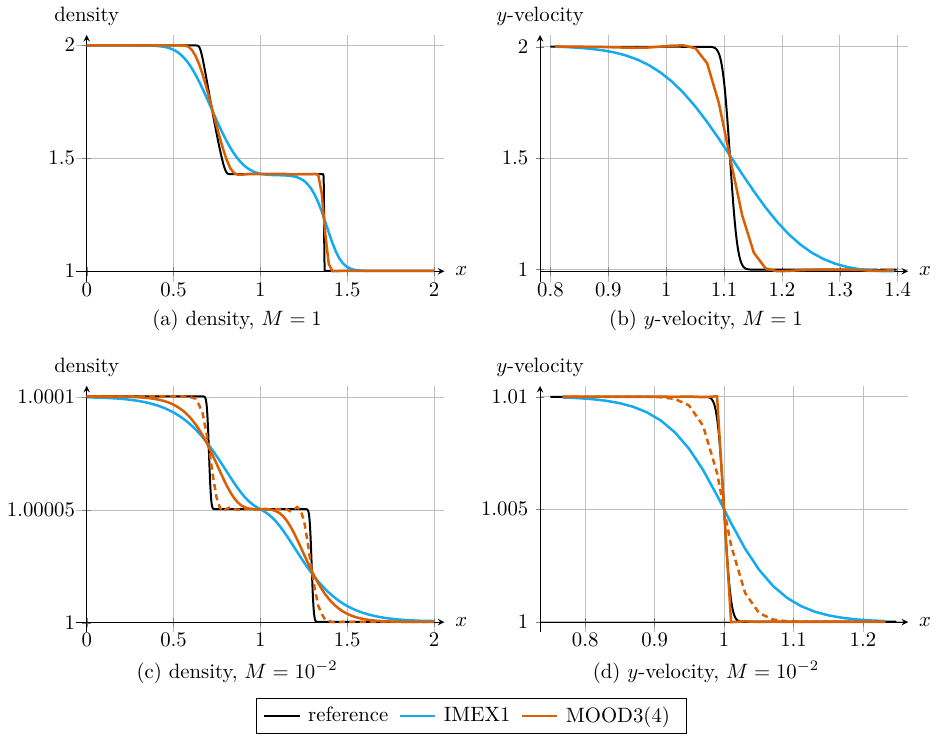}

	}

	\caption{
		Shear wave experiment:
		Approximation of the solution to the Riemann problem~\eqref{eq:Euler_shear_wave}
		at time $T_f = 0.25 M$ on a $100 \times 3$ mesh,
		with $M = 1$ (top)
		and $M = 10^{-2}$ (bottom),
		using the IMEX$1$ and MOOD$3$($4$) schemes.
		For $M = 10^{-2}$, we report the results with material $\nu_{\text{mat}} = 0.1$ (solid lines) and acoustic $\nu_{\text{ac}} = 0.5$ (dashed lines) CFL restrictions.
	}

	\label{fig:Euler_material_wave}

\end{figure}

Since the Riemann problem in the previous test consisted only of acoustic waves, we now consider the approximation of a material wave by introducing a shear wave triggered by a non-zero initial velocity in $y$-direction.
The initial condition is given by
\begin{equation}
	\label{eq:Euler_shear_wave}
	\rho(x,y,0) =
	\begin{dcases}
		1 + M^2 & \text{ if } x < 1, \\
		1       & \text{ otherwise,}
	\end{dcases},
	\quad
	\bm u(x,y,0) =
	\begin{dcases}
		\begin{pmatrix} 0 \\ 1 + M \end{pmatrix} & \text{ if } x < 1, \\
		\begin{pmatrix} 0 \\ 1 \end{pmatrix}     & \text{ otherwise.}
	\end{dcases}
\end{equation}
This experiment verifies that our scheme is able to provide a sharp approximation of the slow-moving shear wave.
The numerical results are given in Figure~\ref{fig:Euler_material_wave} for $M =1$ and $M = 10^{-2}$,
with $\nu_{\text{ac}} = 0.5$ and $\nu_{\text{mat}} = 0.1$,
respectively corresponding to $\Delta t_{\text{ac}} \simeq 8.38 \times 10^{-5}$ and $\Delta t_{\text{mat}} = 9.90 \times 10^{-4}$ for $M = 10^{-2}$,
at the final time $T_f = 0.25 M$.

We observe, as expected, that the approximation of the shear wave is sharp, especially for the material CFL condition.
The acoustic waves in the density solution are diffused for the material CFL, which is expected.
For the acoustic CFL, all waves are captured with the correct propagation speed.
The MOOD procedure is activated respectively $62\%$ and $24\%$ of time iterations for both Mach number regimes using an acoustic CFL condition.
Here, the violation of the detection criterion \eqref{eq:MOODcrit} leads to the use of the lowest stage of the MOOD procedure, which is based on the TVD time integration and TVD space discretisation.


		\paragraph{Circular explosion problem}
		\label{sec:Euler_explosion}
		To demonstrate the performance of the MOOD scheme regarding the capturing of shock waves in the compressible regime, we consider a circular explosion problem.
		The initial condition was originally given in the context of a circular dam-break for the shallow water equations in \cite{Lev2002}.
		To analyse the symmetry of the numerical solution, we use a radial set-up on the domain $\Omega = (-0.5, 0.5)^2$ with $\gamma = 1.4$ and $M=1$.
		The initial data for a radius $r = \sqrt{x^2 + y^2}$ is given by
		\begin{equation}
			\label{eq:Explosion_init}
			\rho(r,0) = \begin{cases}
				2 & \text{ for } r < 0.2 \\
				1 & \text{ otherwise.}   \\
			\end{cases},
			\quad u(r,0) = 0.
		\end{equation}
		The simulation is performed on a uniform grid with $100$ cells in each space direction up to a final time $T_f = 0.125$.
		For simplicity, we have applied periodic boundary conditions since the shock front is contained in the computational domain during the whole simulation.

\begin{figure}[!ht]

	\centering
	\makebox[1.\linewidth][c]{

		\includegraphics[height=6cm,width=\textwidth]{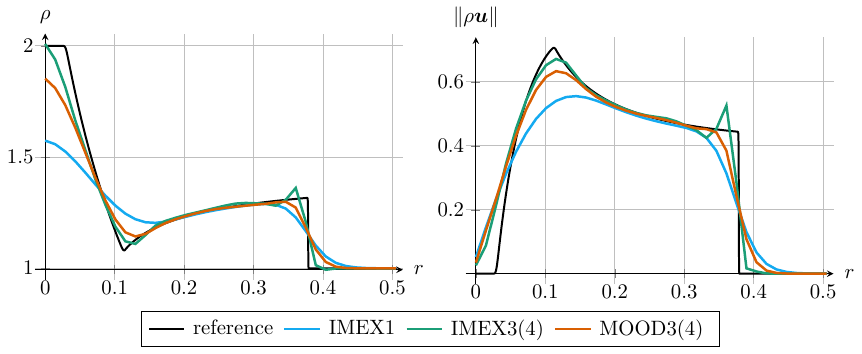}

	}

	\caption{%
		Results for the circular explosion problem from Section~\ref{sec:Euler_explosion},
		depicted along the diagonal $y = x$ line, for $M=1$.
		The radius $r$ corresponds to the position along this diagonal line:
		$r = 0$ corresponds to the point $(0, 0)$,
		while $r = 0.5$ corresponds to the point $(\sqrt{2}/4, \sqrt{2}/4)$.
		Left panel: density $\rho$.
		Right panel: momentum norm $\lVert \rho \boldsymbol u \rVert$.
	}

	\label{fig:Euler_explosion_1D}

\end{figure}

\begin{figure}[!ht]

	\centering
	\makebox[0.95\textwidth][c]{

		\includegraphics{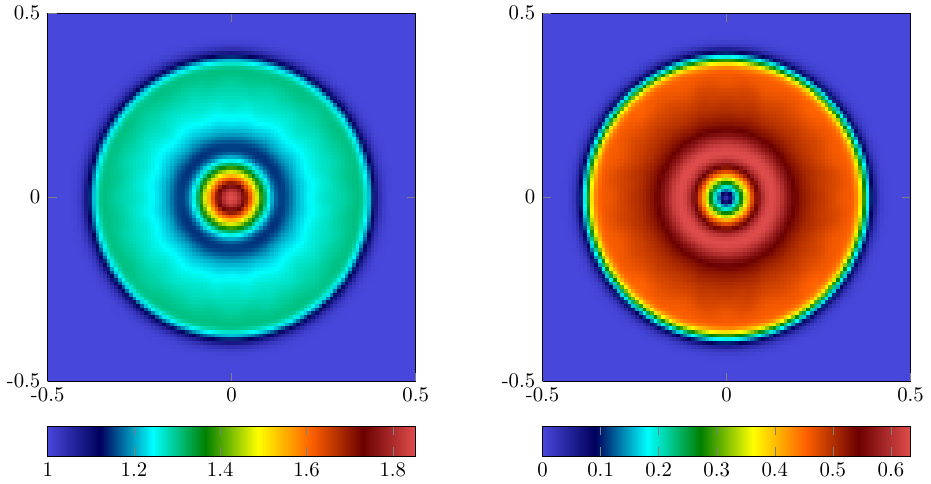}

	}

	\caption{%
		Results for the circular explosion problem from Section~\ref{sec:Euler_explosion},
		depicted in the $(x, y)$ plane for $M=1$.
		Left panel: density $\rho$.
		Right panel: momentum norm $\lVert \rho \boldsymbol u \rVert$.
	}

	\label{fig:Euler_explosion_2D}

\end{figure}

		From the numerical results,
		which are given in Figure~\ref{fig:Euler_explosion_2D},
		we see that the solution is symmetric despite using a Cartesian mesh,
		and that the shock front is captured accurately.
		To see the advantage of a MOOD scheme, we compare the first order IMEX$1$, the high order IMEX$3$($4$) and our MOOD$3$($4$) scheme with
		a reference solution obtained by solving the radial 1D Riemann Problem \eqref{eq:Explosion_init}
		using a second-order minmod HLL scheme with 5000 points following the model in~\cite{CriBal2018}.
		The results are displayed in Figure \ref{fig:Euler_explosion_1D}.
		We observe, as expected, that the IMEX$1$ scheme is quite diffusive and the high order IMEX$3$($4$) overshoots around the shock wave travelling to the right,
		and even in smooth regions between the rarefaction wave and the shock wave.
		This overshoot and the slight oscillations around the reference solution can be avoided with the MOOD$3$($4$) scheme,
		as the TVD scheme is activated when an oscillation is detected.
		For a total of $24$ time steps, the MOOD procedure was activated $9$ times, of which $8$ required the parachute scheme,
		mostly towards the end of the simulation.
		Moreover, the part of the solution depicted in Figure~\ref{fig:Euler_explosion_1D}
		is cut out of the solution on the whole grid given in Figure \ref{fig:Euler_explosion_2D} and consists of only $36 = \lceil 50/\sqrt{2} \rceil$ cells.
		This demonstrates that our MOOD scheme yields good results even on coarse grids.

\subsubsection{Double shear layer}
\label{sec:Euler_double_shear}

The next test case concerns the validation of the asymptotic preserving property in two dimensions, by simulating a double shear layer.
This test case allows a direct comparison with the incompressible flow solution.
The set-up is taken from \cite{DimLouMicVig2018} and initially described in \cite{WeiShu1994}.
The initial data is well-prepared in the sense of \cite{Del2010}, since the initial density is constant and the velocity field is divergence-free.
They are given by
\begin{equation}
	\rho(x,y,0)=\frac{\pi}{15}, \quad
	\bm u(x,y,0) =
	\begin{cases}
		\begin{pmatrix}
			\tanh((y-\pi/2)/\rho) \\
			0.05 \sin(x)
		\end{pmatrix}    & \text{if } y \leq \pi, \\[10pt]
		\begin{pmatrix}
			\tanh((3\pi/2 - y)/\rho) \\
			0.05 \sin(x)
		\end{pmatrix} & \text{otherwise.}
	\end{cases}
\end{equation}
The computational domain is periodic with $\Omega = (0, 2\pi)^2$ and we set $\gamma = 1.4$.

To compute the reference solution given by the single Mach number limit, we follow the procedure in \cite{DimLouMicVig2018} to solve the incompressible Euler equations on a uniform $500 \times 500$ grid with first-order accuracy.
We run the simulation up to a final time of $T_f = 10$ for several Mach numbers $M=10^{-k}$, $k = 1, \ldots, 7$, to validate that our scheme which is based on approximating compressible equations preserves the asymptotic limit towards the incompressible Euler equations.
The results are computed with the MOOD$3$($4$) scheme on a uniform $200 \times 200$ grid, with grid size $\Delta x = 2 \pi / 150 \simeq 3.14 \cdot 10^{-2}$.
Since the problem is sufficiently smooth, the very diffusive parachute scheme is not activated.

As we focus with this test case on the convergence with respect to the Mach number towards an incompressible limit solution and not the mesh convergence, we have computed in Table~\ref{tab:Euler_double_shear_layer_Density_wrt_Mach} the density error against the averaged density $\bar \rho(T_f) = \int_\Omega \rho(x, y, T_f) d (x,y)$.
We see that the density converges in accordance with well-prepared data with an order of $M^2$ against the constant incompressible limit density.
Note that for $M=10^{-7}$ the error is basically of machine precision.

\begin{figure}[!ht]

	\centering
	\makebox[0.95\textwidth][c]{

		\includegraphics{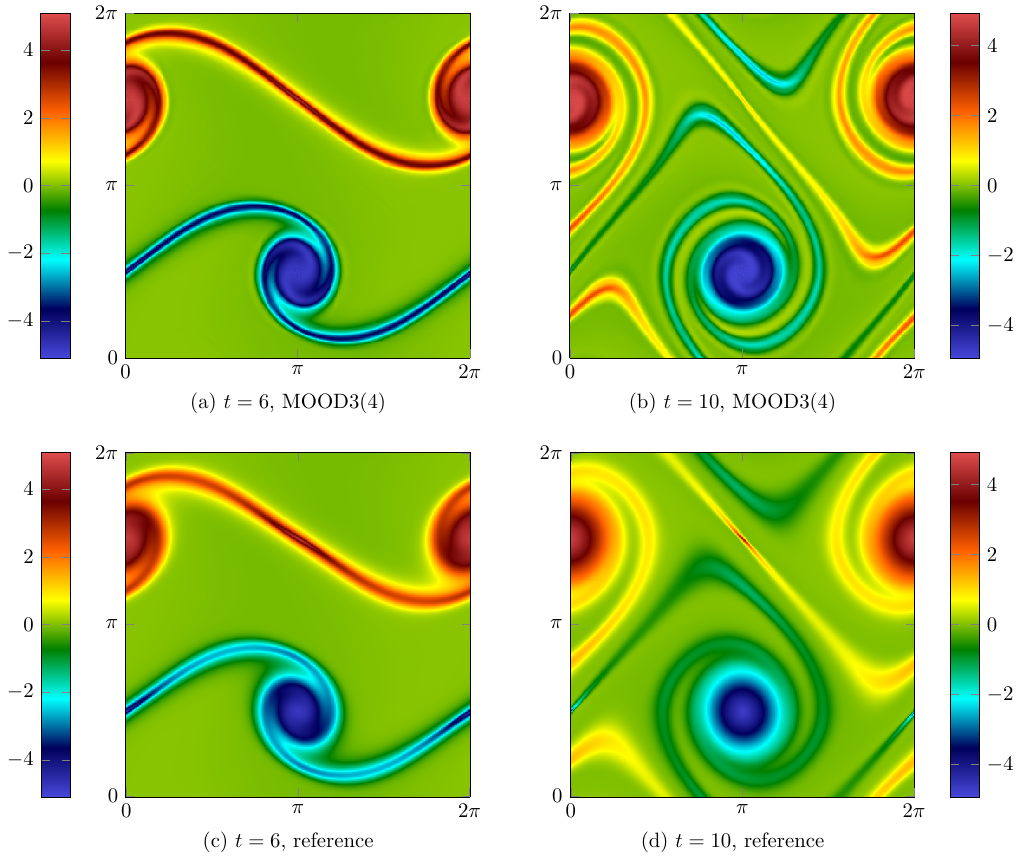}

	}

	\caption{%
		Vorticity $\omega$ for the double shear layer from Section~\ref{sec:Euler_double_shear},
		depicted in the $(x, y)$ plane.
		Left panels: results at $t = 6$.
		Right panels: results at $t = 10$.
		Top panels: MOOD$3$($4$) scheme with $M = 10^{-6}$.
		Bottom panels: reference solution.
	}

	\label{fig:Euler_double_shear_layer}

\end{figure}

In Figure \ref{fig:Euler_double_shear_layer}, the vorticity for the incompressible reference solution and nearly incompressible flow in a Mach number regime of $M=10^{-6}$ is displayed for $T_f = 6$ and $T_f = 10$.
The solution obtained with the compressible isentropic Euler equations shows good agreement with the reference solution and sharply displays the vorticity of the flow.
Note that the incompressible solver is only first order accurate and thus more diffusive than the third order scheme in the MOOD procedure.

\begin{table}[!ht]
		\centering
		\begin{tabular}{cccccccc}
			\toprule
			$M$         & $10^{-1}$            & $10^{-2}$            & $10^{-3}$            & $10^{-4}$            & $10^{-5}$             & $10^{-6}$             & $10^{-7}$             \\
			\cmidrule(rl){1-8}
			$L^2$ error & $2.94 \cdot 10^{-2}$ & $2.98 \cdot 10^{-4}$ & $2.96 \cdot 10^{-6}$ & $2.97 \cdot 10^{-8}$ & $2.97 \cdot 10^{-10}$ & $2.97 \cdot 10^{-12}$ & $2.94 \cdot 10^{-14}$ \\
			\bottomrule
		\end{tabular}

		\caption{%
			$L^2$ errors between density and averaged density
			for the double shear layer from Section~\ref{sec:Euler_double_shear},
			with respect to the Mach number~$M$,
			at time $T_f = 10$,
			for a mesh with $25 \times 25$ cells.
		}

		\label{tab:Euler_double_shear_layer_Density_wrt_Mach}
\end{table}

\subsubsection{Stationary isentropic vortex}
\label{sec:Euler_vortex}

We finally consider a stationary two-dimensional vortex, given by
\begin{equation}
	\label{eq:vortex}
	\begin{dcases}
		\rho(x, y) = 1 - \frac{M^2}{8} e^{- 2 a^2 r(x, y)^2}, \\
		\bm u(x,y) = a \sqrt{\frac \gamma 2} e^{- a^2 r(x, y)^2} \rho(x, y)^{\frac{\gamma}{2} - 1} \begin{pmatrix} y \\ -x \end{pmatrix}.
	\end{dcases}
\end{equation}
For the simulation, we set $a = 8$ and consider the computational domain $(0, 1)^2$ with periodic boundary conditions up to the final time $T_f = 0.2$.
We take $\nu_{\text{mat}} = 0.1$, which corresponds to $\nu_{\text{ac}} \simeq 0.124$ for $M = 1$ and $\nu_{\text{ac}} \simeq 19.2$ for $M = 10^{-2}$.
The time step is then given by $\Delta t = 0.01 \sqrt{N}/32$.

\begin{figure}[!ht]

	\centering
	\makebox[0.95\linewidth][c]{

		\includegraphics{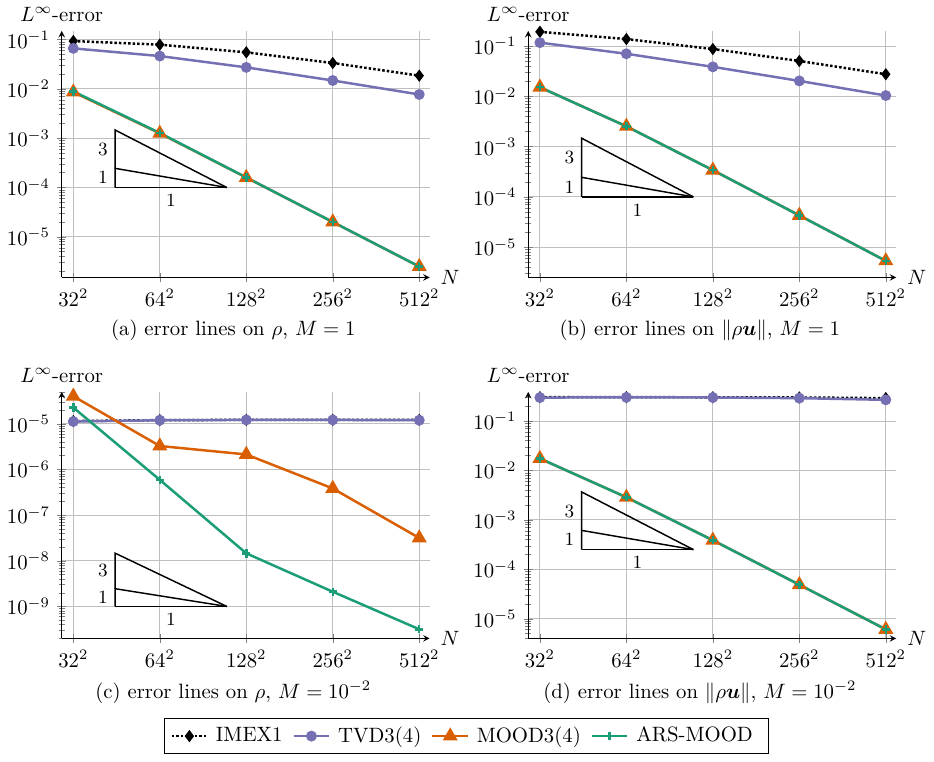}

	}

	\caption{
		Error lines in $L^2$ norm for the 2D vortex described in Section~\ref{sec:Euler_vortex},
		with $M = 1$ (top panels) and $M = 10^{-2}$ (bottom panels).
		Left panels: errors on the density~$\rho$;
		right panels: errors on the momentum norm $\lVert \rho \boldsymbol u \rVert$.
	}

	\label{fig:Euler_vortex_error_lines}

\end{figure}

The $L^\infty$ error lines for $M = 1$ and $M = 10^{-2}$
are depicted in Figure~\ref{fig:Euler_vortex_error_lines},
		and the numerical values of the errors and of the orders of accuracy
		are reported in
		Table~\ref{tab:vortex_errors_and_orders}.%
We note that the TVD3 scheme is first-order accurate and more precise than the IMEX1 scheme for $M = 1$.
For $M=10^{-2}$, the TVD3(4) and IMEX1 schemes have roughly the same error.
This is due to the Mach number-dependent diffusion required required for TVD stability in the numerical scheme, see \cite{DimLouVig2017}.
This highlights the relevance of our three-stage MOOD algorithm.
Further, we note that the MOOD3(4) scheme has the expected EOC of third order.
Since the solution is smooth, a reduction of the EOC due to the use of the parachute scheme has not occurred.
This confirms the fact that the detection criterion based on the Riemann invariants is able to correctly detect smooth regions, and to appropriately keep the high-order schemes when the solution is oscillation-free.
As a consequence, the TVD-MOOD scheme remains third-order accurate for each value of the Mach number $M$ under consideration.
Furthermore, we see from the errors that the diffusion is independent of the Mach number, and since the initial condition is well-prepared, see \cite{Del2010}, it is a strong indication that the scheme is also asymptotic preserving, i.e. is also a consistent discretization of the incompressible Euler equations as $M$ tends to 0.

We also report the error lines of the ARS-MOOD scheme,
which consists of performing the MOOD Algorithm~\ref{algo:implicit_MOOD}
with the ARS(2,3,3) as the third order scheme and the TVD3(4) time integrator for the parachute scheme.
We note that the ARS-MOOD performs as well as the MOOD3(4) scheme for the errors in the momentum,
while it yields significantly better density errors for $M = 10^{-2}$.
This can be explained by the increased stability properties of the ARS(2,3,3) leading to a better approximation for small Mach numbers.
Although both MOOD3(4) and ARS-MOOD show a third-order EOC, the improved errors of the ARS-MOOD scheme underline the importance of choosing a suitable $L$-stable or stiffly accurate high order scheme in the MOOD procedure for simulating low Mach flows.

\pgfplotstableset{
	points column/.style={%
			/pgfplots/table/display columns/#1/.style={%
					int detect, column name=$N$,
				}
		},
	error column/.style={%
			/pgfplots/table/display columns/#1/.style={%
					sci, sci zerofill, precision=2, sci 10e, column name={$L^2$ error}
				}
		},
	order column/.style={%
			/pgfplots/table/display columns/#1/.style={%
					fixed, fixed zerofill, precision=2, column name={EOC},
				}
		}
}

\begin{table}[!ht]
	\begin{subtable}{\textwidth}\centering
		{
			\pgfplotstabletypeset[
				clear infinite,
				empty cells with={---},
				points column/.list={0},
				error column/.list={1,3,5,7},
				order column/.list={2,4,6,8},
				every head row/.style={
						before row={
								\toprule%
								& \multicolumn{2}{c}{$\rho$, IMEX$1$}
								& \multicolumn{2}{c}{$\rho$, TVD$3$($4$)}
								& \multicolumn{2}{c}{$\rho$, MOOD$3$($4$)}
								& \multicolumn{2}{c}{$\rho$, ARS-MOOD} \\
								\cmidrule(lr){2-3}\cmidrule(lr){4-5}\cmidrule(lr){6-7}\cmidrule(lr){8-9}
							}
					},
				every head row/.append style={
						after row={
								\cmidrule(lr){1-1}\cmidrule(lr){2-3}\cmidrule(lr){4-5}\cmidrule(lr){6-7}\cmidrule(lr){8-9}
							}
					},
				every last row/.style={
						after row=\bottomrule
					},
			]{
				1024   0.09334132 nan        0.06629381 nan        0.00861349     nan        0.00889327     nan
				4096   0.07893188 0.24190776 0.04633344 0.51682037 0.00125513     2.77876111 0.00127633     2.80071259
				16384  0.05515124 0.51721479 0.02730796 0.76273225 0.00015960     2.97530424 0.00016006     2.99531667
				65536  0.03346743 0.72063555 0.01480695 0.88304705 1.99832393e-05 2.99759828 2.00092778e-05 2.99987181
				262144 0.01853669 0.85237412 0.00770678 0.94207439 2.49839446e-06 2.99971728 2.50002315e-06 3.00065574
			}
		}
		\caption{%
			Error of the density $\rho$ for $M=1$.
		}
		\medskip
	\end{subtable}
	\begin{subtable}{\textwidth}\centering
		{
			\pgfplotstabletypeset[
				sci zerofill,
				clear infinite,
				empty cells with={---},
				points column/.list={0},
				error column/.list={1,3,5,7},
				order column/.list={2,4,6,8},
				every head row/.style={
						before row={
								\toprule%
								& \multicolumn{2}{c}{$\lVert \rho \boldsymbol u \rVert$, IMEX$1$}
								& \multicolumn{2}{c}{$\lVert \rho \boldsymbol u \rVert$, TVD$3$($4$)}
								& \multicolumn{2}{c}{$\lVert \rho \boldsymbol u \rVert$, MOOD$3$($4$)}
								& \multicolumn{2}{c}{$\lVert \rho \boldsymbol u \rVert$, ARS-MOOD} \\
								\cmidrule(lr){2-3}\cmidrule(lr){4-5}\cmidrule(lr){6-7}\cmidrule(lr){8-9}
							}
					},
				every head row/.append style={
						after row={
								\cmidrule(lr){1-1}\cmidrule(lr){2-3}\cmidrule(lr){4-5}\cmidrule(lr){6-7}\cmidrule(lr){8-9}
							}
					},
				every last row/.style={
						after row=\bottomrule
					},
			]{
				1024   0.19197562 nan        0.11685122 nan        0.01517031     nan        0.01515047     nan
				4096   0.13749262 0.48156892 0.07016972 0.73575229 0.00253258     2.58257082 0.00253033     2.581965097966
				16384  0.08717320 0.65739761 0.03858335 0.86287019 0.00034026     2.89589837 0.00034002     2.89563404
				65536  0.05045620 0.78885311 0.02024052 0.93073206 4.28632057e-05 2.98882591 4.28299531e-05 2.98892761
				262144 0.02748182 0.87655401 0.01036921 0.96494037 5.35140559e-06 3.00174997 5.34746542e-06 3.00169294
			}
		}
		\caption{%
			Error of the momentum norm $\lVert \rho \boldsymbol u \rVert$ for $M=1$.
		}
		\medskip
	\end{subtable}
	\begin{subtable}{\textwidth}\centering
		{
			\pgfplotstabletypeset[
				sci zerofill,
				clear infinite,
				empty cells with={---},
				points column/.list={0},
				error column/.list={1,3,5,7},
				order column/.list={2,4,6,8},
				every head row/.style={
						before row={
								\toprule%
								& \multicolumn{2}{c}{$\rho$, IMEX$1$}
								& \multicolumn{2}{c}{$\rho$, TVD$3$($4$)}
								& \multicolumn{2}{c}{$\rho$, MOOD$3$($4$)}
								& \multicolumn{2}{c}{$\rho$, ARS-MOOD} \\
								\cmidrule(lr){2-3}\cmidrule(lr){4-5}\cmidrule(lr){6-7}\cmidrule(lr){8-9}
							}
					},
				every head row/.append style={
						after row={
								\cmidrule(lr){1-1}\cmidrule(lr){2-3}\cmidrule(lr){4-5}\cmidrule(lr){6-7}\cmidrule(lr){8-9}
							}
					},
				every last row/.style={
						after row=\bottomrule
					},
			]{
				1024   1.14358674e-05 nan         1.12087401e-05  nan        3.99854036e-05 nan        2.26277657e-05 nan
				4096   1.19994140e-05 -0.06939816 1.19362408e-05 -0.09072442 3.25904996e-06 3.61695008 5.90494304e-07 5.26002718
				16384  1.21442956e-05 -0.01731486 1.21230271e-05 -0.02240144 2.12778637e-06 0.61509816 1.45477014e-08 5.34306002
				65536  1.21759104e-05 -0.00375084 1.21282482e-05 -0.00062120 3.83072547e-07 2.47366377 2.10141771e-09 2.79135635
				262144 1.21416047e-05  0.00407054 1.19216796e-05  0.02478368 3.15889971e-08 3.60012360 3.16837667e-10 2.72954720
			}
		}
		\caption{%
		Error of the density $\rho$ for $M=10^{-2}$.
		}
		\medskip
	\end{subtable}
	\begin{subtable}{\textwidth}\centering
		{
			\pgfplotstabletypeset[
				sci zerofill,
				clear infinite,
				empty cells with={---},
				points column/.list={0},
				error column/.list={1,3,5,7},
				order column/.list={2,4,6,8},
				every head row/.style={
						before row={
								\toprule%
								& \multicolumn{2}{c}{$\lVert \rho \boldsymbol u \rVert$, IMEX$1$}
								& \multicolumn{2}{c}{$\lVert \rho \boldsymbol u \rVert$, TVD$3$($4$)}
								& \multicolumn{2}{c}{$\lVert \rho \boldsymbol u \rVert$, MOOD$3$($4$)}
								& \multicolumn{2}{c}{$\lVert \rho \boldsymbol u \rVert$, ARS-MOOD} \\
								\cmidrule(lr){2-3}\cmidrule(lr){4-5}\cmidrule(lr){6-7}\cmidrule(lr){8-9}
							}
					},
				every head row/.append style={
						after row={
								\cmidrule(lr){1-1}\cmidrule(lr){2-3}\cmidrule(lr){4-5}\cmidrule(lr){6-7}\cmidrule(lr){8-9}
							}
					},
				every last row/.style={
						after row=\bottomrule
					},
			]{
				1024   3.02980146e-01  nan        2.98326013e-01  nan         1.74027668e-02  nan        1.74284028e-02 nan
				4096   3.03212991e-01  -0.00110831 3.01406996e-01  -0.01482312  2.86335503e-03 2.60353823 2.86331977e-03 2.60567966
				16384  3.02618765e-01 0.00283012 2.99194023e-01 0.01063155  3.87340364e-04 2.88603281 3.86634380e-04 2.88864696
				65536  2.99562167e-01 0.01464602 2.90452919e-01 0.04277703  4.87871580e-05 2.98902850 4.87411543e-05 2.98775761
				262144 2.90573275e-01 0.04395341 2.68352756e-01 0.11417362  6.08341368e-06 3.00354843 6.07906001e-06 3.00322025
			}
		}
		\caption{%
		Error of the momentum norm $\lVert \rho \boldsymbol u \rVert$ for $M=10^{-2}$.
		}
	\end{subtable}
	\caption{%
		Errors in $L^2$ norm for the density $\rho$ and the momentum norm $\lVert \rho \boldsymbol u \rVert$ associated to the error lines given in Figure \ref{fig:Euler_vortex_error_lines} for the 2D vortex described in Section~\ref{sec:Euler_vortex}.
	}
	\label{tab:vortex_errors_and_orders}
\end{table}


\section{Conclusions and future work}
\label{sec:conclusion}
We have presented a new approach on constructing first order TVD IMEX-RK schemes, which are especially suited as MOOD parachute schemes to simulate multi-scale equations due to a reduced numerical viscosity compared to a backward forward first order Euler integrator, see Table~\ref{tab:space_time_errors} for a quantification using space-time errors.
Their development is motivated by the fact that there is a first order barrier for IMEX-TVD schemes that have a scale independent CFL restriction.
The key in constructing those schemes lies in using a convex combination of a first-order TVD IMEX scheme with a high-order IMEX RK scheme.
We gave a theoretical justification of our TVD approach by means of studying a one dimensional linear scalar equation.
The obtained TVD scheme is then used as a parachute scheme in a MOOD procedure, which is triggered whenever the TVD property is violated by the high order scheme.
The performance of the resulting TVD and MOOD schemes was verified by approximating a discontinuous solution for scalar linear transport, see Figure~\ref{fig:step_order_3_solutions}, as well as Riemann problems for the 2D isentropic Euler equations, see Figures~\ref{fig:Euler_RP},~\ref{fig:Euler_material_wave} and \ref{fig:Euler_explosion_1D} assessing the resolution of acoustic and material waves.
		In addition, the scheme is shown to preserve the correct asymptotics in the singular low Mach number limit, see Figure~\ref{fig:Euler_double_shear_layer} and Table~\ref{tab:Euler_double_shear_layer_Density_wrt_Mach}.
Further, it was verified that the MOOD procedure yields high order convergence, correctly detecting smooth solutions, by simulating a vortex in different Mach number regimes, see Figure~\ref{fig:Euler_vortex_error_lines}.
Our results are a significant improvement to the scheme from~\cite{DimLouMicVig2018} for the isentropic Euler equations, especially for small Mach numbers, where our MOOD schemes are able keep a sharp profile on the material wave.
Due to the Mach number dependent diffusion of the TVD scheme from \cite{DimLouMicVig2018}, the grid would have to be drastically refined to obtain comparable results on the material wave.

The construction of TVD or SSP IMEX schemes with a material CFL restriction for multi-scale equations is still an active field of research, as shown by the recent work of \cite{GotGraHuShu2022}.
Therefore, our schemes, which provide a certain flexibility regarding the time step while maintaining the TVD property, are necessary to obtain physically admissible solutions for problems relevant to the community.

Since, in this work, we have neglected a higher order reconstruction of the implicitly treated derivatives due to avoiding the inversion of non-linear systems, in future work we plan to combine our schemes with the linear implicit high order \textit{Quinpi} approach from \cite{PupSemVis2022}.


\appendix

\section{On non-CK IMEX schemes}
\label{sec:appendix_non_CK_schemes}

Consider the following Butcher tableaux, defining an IMEX scheme in non-CK, non-ARS form:
\begin{equation}
	\label{tab:Butcher_tableaux_non_CK}
	\text{explicit: }
	\renewcommand{\arraystretch}{1.25}
	\begin{array}{c|ccccc}
		0          & 0             & 0      & \cdots           & 0            \\
		\tilde c_2 & \tilde a_{21} & 0      & \cdots           & 0            \\
		\vdots     & \vdots        & \ddots & \ddots           & \vdots       \\
		\tilde c_s & \tilde a_{s1} & \cdots & \tilde a_{s,s-1} & 0            \\ \hline
		           & \tilde b_1    & \cdots & \tilde b_{s-1}   & \tilde b_{s}
	\end{array}
	\hskip1cm
	\text{implicit: }
	\renewcommand{\arraystretch}{1.25}
	\begin{array}{c|cccc}
		c_1    & a_{11} & 0      & \cdots & 0      \\
		c_2    & a_{21} & a_{22} & \cdots & 0      \\
		\vdots & \vdots & \vdots & \ddots & \vdots \\
		c_s    & a_{s1} & a_{s2} & \cdots & a_{ss} \\ \hline
		       & b_{1}  & b_{2}  & \cdots & b_{s}
	\end{array}.
\end{equation}
We derive stability conditions analogous to Theorem~\ref{theo:TVDconvex} for this case where the first column of the implicit tableau is non-zero.
After lengthy computations, we get the following result:
\begin{theorem}
	\label{theo:TVDconvex_nonCK}
	Let $\tilde A, A \in \mathbb{R}^{s\times s}$, $\tilde b, b, \tilde c, c \in \mathbb{R}^s$
	define two Butcher tableaux~\eqref{tab:Butcher_tableaux_non_CK} fulfilling~\eqref{cond:c} and the $p$-th order compatibility conditions.
	Let $\tilde b$ and $b$ coincide with the last rows of $\tilde A$ and $A$ respectively.
	For $k = 1, \ldots, s$ and $l = 1, \ldots, k-1$, we define
	\begin{equation*}
		\mathcal{A}_k = \theta_k a_{kk} + (1 - \theta_k)c_k, \quad
		\mathcal{\tilde A}_k = (1-\theta_k) \tilde c_k, \quad
		\mathcal{B}_{kl} = \frac{\theta_k a_{kl}}{\mathcal{A}_l}, \quad
		\mathcal{\tilde B}_{kl} = \theta_k \tilde a_{kl}.
	\end{equation*}
	In addition, we recursively define the following expressions:
	\begin{alignat*}{3}
		 & \mathcal{\tilde C}_k = \mathcal{\tilde A}_k - \sum_{l=2}^{k-1} \mathcal{B}_{kl} \mathcal{\tilde C}_l, &                                                                                                                  & \qquad
		 &                                                                                                       & \mathcal{\tilde D}_{kl} = \mathcal{\tilde B}_{kl} - \sum_{r=l+1}^{k-1} \mathcal{B}_{kr} \mathcal{\tilde D}_{rl},           \\
		 & \mathcal{C}_k = 1 - \sum_{l=1}^{k-1} \mathcal{B}_{kl} \mathcal{C}_l,                                  &                                                                                                                  & \qquad
		 &                                                                                                       & \mathcal{D}_{kl} = \mathcal{B}_{kl} - \sum_{r=l+1}^{k-1} \mathcal{B}_{kr} \mathcal{D}_{rl}.
	\end{alignat*}
	Then, under the following restrictions for $k = 1, \ldots, s$ and $l = 1, \ldots, k-1$,
	\begin{equation*}
		\mathcal{A}_k > 0,
		\qquad
		0 \leq \lambda \mathcal{\tilde C}_k \leq \mathcal{C}_k,
		\qquad
		0 \leq \lambda \mathcal{\tilde D}_{k,l} \leq \mathcal{D}_{k,l},
	\end{equation*}
	the scheme consisting in the convex combination based on the Butcher tableaux~\eqref{tab:Butcher_tableaux_non_CK},
	combined with a TVD limiter, is $L^\infty$ stable and TVD under a CFL condition determined by $\lambda \geq 0$ where $\lambda$ does not depend on $\varepsilon$.
\end{theorem}

When performing numerical experiments, we observe that the results of schemes derived under the conditions of Theorem~\ref{theo:TVDconvex_nonCK}
are not as compelling as results of schemes obeying Theorem~\ref{theo:TVDconvex}.
Therefore, we do not include such schemes in the numerical experiments,
but we still state Theorem~\ref{theo:TVDconvex_nonCK} for the sake of completeness.

\section{
  \texorpdfstring{TVD$3$($4$)}{TVD3(4)}
 }
\label{sec:appendix_TVD34}

For the TVD$3$($4$) scheme, the explicit Butcher tableau is given by
\begin{equation*}
	\label{eq:TVD3_4_explicit_tableau}
	\renewcommand{\arraystretch}{1.25}
	\resizebox{\linewidth}{!}{%
		$
			\begin{array}{c|cccc}
				0                  & 0                   & 0                  & 0                  & 0                   \\
				0.2049503677289891 & 0.2049503677289891  & 0                  & 0                  & 0                   \\
				0.4173127343286904 & 0.2123925641886599  & 0.2049201701400305 & 0                  & 0                   \\
				0.9048203025659662 & -0.4501877125339555 & 0.3955748607480934 & 0.9594331543518283 & 0                   \\
				\hline
				                   & 0                   & 0.3354718384287510 & 0.3487815573407456 & 0.3157466042305059
			\end{array}
		$}
\end{equation*}
while the implicit Butcher tableau is given by
\begin{equation*}
	\label{eq:TVD3_4_implicit_tableau}
	\renewcommand{\arraystretch}{1.25}
	\begin{array}{c|cccc}
		0                  & 0 & 0                  & 0                  & 0                   \\
		0.2049503677289891 & 0 & 0.2049503677289891 & 0                  & 0                   \\
		0.4173127343286904 & 0 & 0.2040104873103189 & 0.2133022470183705 & 0                   \\
		0.9048203025659662 & 0 & 0.3991926529002874 & 0.4115004113464103 & 0.0941272383192684  \\
		\hline
		                   & 0 & 0.3354718384287510 & 0.3487815573407456 & 0.3157466042305059
	\end{array}.
\end{equation*}
with the parameters
\begin{align*}
	 & \theta_1 = 1, \quad \theta_2 = 1, \quad \theta_3 = 1, \quad \theta_4 = 0.5110907014643069, \quad \theta_5 = 0.4997722865197203, \\
	 & \lambda = 0.5471076190680170.
\end{align*}

\section{
  \texorpdfstring{ARS($2,2,3$)}{ARS(2,2,3)}
 }
\label{sec:appendix_ARS223}

The ARS(2,2,3) scheme from~\cite{AscRuuSpi1997} is given by the following tableaux:
\begin{equation*}
	\begin{array}{c|ccc}
		0          & 0          & 0            & 0         \\
		\delta     & \delta     & 0            & 0         \\
		1 - \delta & \delta - 1 & 2 - 2 \delta & 0         \\
		\hline
		           & 0          & \sfr 1 2     & \sfr 1 2

	\end{array},
	\quad
	\begin{array}{c|ccc}
		0          & 0 & 0            & 0         \\
		\delta     & 0 & \delta       & 0         \\
		1 - \delta & 0 & 1 - 2 \delta & \delta    \\
		\hline
		           & 0 & \sfr 1 2     & \sfr 1 2

	\end{array},
	\quad \text{with } \delta = \frac{3 + \sqrt{3}}{6}.
\end{equation*}


\section*{Acknowledgements}

This work has been partially funded by a CNRS-INSMI PEPS JCJC project.
A.T. has been partially supported by the Gutenberg Research College, JGU Mainz.


\bibliographystyle{plain}
\bibliography{./literature}

\end{document}